\documentclass{IEEEtran}

\usepackage{enumitem,kantlipsum}
\usepackage{graphicx}
\usepackage{epstopdf}
\usepackage{amssymb}
\usepackage{amsmath}
\usepackage{subfigure}
\usepackage{epsfig}
\usepackage{color}
\usepackage{enumitem}
\usepackage{float}
\usepackage{soul}
\usepackage{wrapfig}
\usepackage{arydshln}
\usepackage{tablefootnote}

\usepackage{amsthm}
\def\0{{\bf 0}}
\def\1{{\bf 1}}

\def\beq{\begin{equation*}}
    \def\eeq{\end{equation*}}
\def\bql{\begin{equation}}
    \def\eql{\end{equation}}
\def\bqn{\begin{eqnarray*}}
    \def\eqn{\end{eqnarray*}}
\def\bnl{\begin{eqnarray}}
    \def\enl{\end{eqnarray}}
\def\bma{\begin{bmatrix}}
    \def\ema{\end{bmatrix}}
\def\bmx{\begin{matrix}}
    \def\emx{\end{matrix}}
\def\ben{\begin{enumerate}}
    \def\een{\end{enumerate}}
\def\bit{\begin{itemize}}
    \def\eit{\end{itemize}}
\def\bei{\begin{itemize}}
    \def\eei{\end{itemize}}
\def\bet{\begin{tabular}}
    \def\eet{\end{tabular}}

\newcommand{\ba}{\mathbf{a}}

\newcommand{\R}{\mathbb{R}}

\newcommand{\be}{\mathbf{e}}

\newcommand{\g}{\mathbf{g}}

\newcommand{\barx}{\overline{x}}

\newcommand{\bu}{\mathbf{u}}

\newcommand{\bv}{\mathbf{v}}

\usepackage[dvipsnames]{xcolor}

\def\R{\mathbb{R}}

\def\1{{\bf1}}
\def\la{\langle}
\def\ra{\rangle}

\def\a{\alpha}
\def\g{\gamma}

\def\bit{\begin{itemize}}
\def\eit{\end{itemize}}

\def\be{\begin{equation}}
\def\ee{\end{equation}}
\def\ba{\begin{eqnarray}}
\def\ea{\end{eqnarray}}

\def\bes{\begin{equation*}}
\def\ees{\end{equation*}}
\def\bas{\begin{eqnarray*}}
\def\eas{\end{eqnarray*}}

\usepackage{url}
\usepackage{verbatim}

\newtheorem{Remark 1}{Remark}
\newtheorem{Remark 2}[Remark 1]{Remark}
\newtheorem{Remark 3}[Remark 1]{Remark}
\newtheorem{Remark 4}[Remark 1]{Remark}
\newtheorem{Remark 5}[Remark 1]{Remark}
\newtheorem{Remark 6}[Remark 1]{Remark}
\newtheorem{Remark 7}[Remark 1]{Remark}

\newtheorem{Lemma 1}{Lemma}
\newtheorem{Lemma 2}[Lemma 1]{Lemma}
\newtheorem{Lemma 3}[Lemma 1]{Lemma}
\newtheorem{Lemma 4}[Lemma 1]{Lemma}
\newtheorem{Lemma 5}[Lemma 1]{Lemma}
\newtheorem{Lemma 6}[Lemma 1]{Lemma}
\newtheorem{Lemma 7}[Lemma 1]{Lemma}

\newtheorem{Assumption 1}{Assumption}
\newtheorem{Assumption 2}[Assumption 1]{Assumption}
\newtheorem{Assumption 3}[Assumption 1]{Assumption}
\newtheorem{Assumption 4}[Assumption 1]{Assumption}
\newtheorem{Definition 1}{Definition}
\newtheorem{Theorem 1}{Theorem}
\newtheorem{Theorem 2}[Theorem 1]{Theorem}
\newtheorem{Theorem 3}[Theorem 1]{Theorem}
\newtheorem{Theorem 4}[Theorem 1]{Theorem}
\newtheorem{Theorem 5}[Theorem 1]{Theorem}
\newtheorem{Theorem 6}[Theorem 1]{Theorem}
\newtheorem{Theorem 7}[Theorem 1]{Theorem}
\newtheorem{Theorem 8}[Theorem 1]{Theorem}
\newtheorem{Theorem 9}[Theorem 1]{Theorem}
\newtheorem{Theorem 10}[Theorem 1]{Theorem}

\title{\LARGE \bf
 Gradient-tracking based Distributed Optimization with Guaranteed Optimality under Noisy Information Sharing
 }

\author{Yongqiang Wang, Tamer Ba{\c{s}}ar% <-this % stops a space
\thanks{ The work of the first author  was supported in part by the National Science Foundation under Grants  ECCS-1912702, CCF-2106293, CCF-2215088, and CNS-2219487.  Research of the second author was supported in part by the ONR MURI Grant N00014-16-1-2710 and in part by the Army Research Laboratory, under Cooperative Agreement Number W911NF-17-2-0196.}
\thanks{Yongqiang Wang is with the Department of Electrical and Computer Engineering, Clemson University, Clemson, SC 29634, USA
{\tt\small{yongqiw}@clemson.edu}
}%
\thanks{Tamer Ba{\c{s}}ar is with the Coordinated Science Lab, University of Illinois
  Urbana-Champaign, Urbana, IL 61801, USA {\tt\small
basar1@illinois.edu}}
  }

\begin{document}

\maketitle
\thispagestyle{empty}
\pagestyle{empty}

%%%%%%%%%%%%%%%%%%%%%%%%%%%%%%%%%%%%%%%%%%%%%%%%%%%%%%%%%%%%%%%%%%%%%%%%%%%%%%%%
\begin{abstract}
Distributed optimization enables networked agents to cooperatively solve a global optimization problem even with each participating agent only having access to a local partial view of the objective function. Despite making significant inroads, most existing results on distributed optimization rely on   noise-free information sharing among the agents, which is problematic when communication channels are noisy, messages are coarsely quantized, or shared information are obscured by additive noise for the purpose of achieving differential privacy. The problem of information-sharing noise is particularly  pronounced in the state-of-the-art gradient-tracking based distributed optimization algorithms, in that   information-sharing noise will accumulate with iterations on the gradient-tracking estimate of these algorithms, and the ensuing variance will even grow  unbounded when the noise is persistent. This paper proposes a new gradient-tracking based distributed optimization approach that can avoid information-sharing noise from accumulating in the gradient estimation. The approach is applicable even when the {inter-agent interaction is} time-varying, which is key to enable the incorporation of a decaying factor in inter-agent interaction to gradually eliminate the influence of information-sharing noise. In fact, we rigorously prove that the proposed approach can
ensure the almost sure  convergence of all agents to the same optimal solution even in the presence of persistent information-sharing noise. The approach is applicable to general directed graphs. It is also capable of ensuring the almost sure convergence of all agents  to an optimal solution when the gradients are noisy, which is common in  machine learning applications. Numerical simulations confirm the effectiveness of the proposed approach.
\end{abstract}

\section{Introduction}
%Main observations:
%\begin{itemize}
%\item In conventional gradient-tracking algorithms, to achieve differential privacy, we have to
%add noise to both $y^k$ and $x^k$, i.e., share $y^k+\xi^k$ and
%$x^k+\zeta^k$;
%\item  The noise will accumulate in gradient estimate, i.e, ${\bf 1}^Ty^k={\bf 1}^T\left(\nabla f(\cdot)+\sum_{\ell=0}^{k-1}\xi^{\ell}\right)$ ;
%
%\end{itemize}
We consider a distributed convex optimization problem where multiple agents cooperatively solve a global optimization problem via local computations and local sharing of information. This is motivated by the problem's broad applications in cooperative control
\cite{yang2019survey}, distributed sensing
\cite{bazerque2009distributed}, multi-agent systems
\cite{raffard2004distributed}, sensor networks
\cite{zhang2017distributed}, and large-scale machine learning
\cite{tsianos2012consensus}.   In many of these applications, each agent  has access to only a local portion of the objective function. Such a distributed optimization problem can be formulated  in the following general form:
 \begin{equation}\label{eq:optimization_formulation1}
\min\limits_{\theta\in\mathbb{R}^d} F(\theta)\triangleq
\frac{1}{m}\sum_{i=1}^m f_i(\theta)
\end{equation}
where $m$ is the number of agents, $\theta\in\mathbb{R}^d$ is
 a decision variable common to all agents, while
$f_i:\mathbb{R}^d\rightarrow\mathbb{R}$ is a local objective
function private to agent $i$.

To solve problem (\ref{eq:optimization_formulation1}) in a distributed manner, plenty of  algorithms have been reported since the seminal works in the 1980s \cite{tsitsiklis1984problems}. Some of the popular algorithms include gradient methods (e.g.,
\cite{nedic2009distributed,shi2015extra,xu2017convergence,qu2017harnessing,xin2018linear}),
distributed alternating direction method of multipliers (e.g.,
\cite{shi2014linear,zhang2019admm}),  and distributed Newton methods
(e.g., \cite{wei2013distributed,zhang2022distributed}). In this paper, we focus on distributed gradient methods, which are particulary appealing for agents with limited computational or storage capabilities due to their low computation complexity and storage requirement. Existing distributed gradient methods can be generally divided into two categories. The first category of distributed gradient  methods  directly concatenate gradient based steps with a consensus operation of the optimization variable (referred to as the static-consensus based approach hereafter), with typical examples including~\cite{nedic2009distributed,yuan2016convergence}. These approaches are simple and efficient in computation since they   require each agent to share only one variable in each iteration. However, they are only applicable in undirected graphs and directed graphs that are  balanced  (the sum of each agent's in-neighbor coupling weights equal to the sum of its out-neighbor coupling weights). The second category of distributed gradient methods   exploit  a dynamic-consensus mechanism  to track the global gradient (and hence usually called gradient-tracking based approach), and are  applicable to general directed graphs (see, e.g., \cite{xu2017convergence,qu2017harnessing,xin2018linear,di2016next,pu2020push}). Such approaches can ensure convergence to an optimal solution under constant  stepsizes and, hence, can achieve faster convergence. However, these approaches need every agent to maintain and share an additional gradient-tracking variable besides the optimization variable, which doubles the communication overhead  compared with the approaches in the first category.

Although plenty of inroads have been made in distributed optimization, most of the existing approaches assume noise-free information sharing, i.e., every agent is able to  acquire neighboring agents' intermediate local optimization variables accurately without any  distortion or corruption. Such an assumption does not hold any longer, however, in many application scenarios. For example, when  communication channels are noisy, every message will be distorted by channel noise which is usually modeled as additive Gaussian noise \cite{kar2008distributed,srivastava2011distributed}. An even more pervasive source of noise in information sharing comes from quantization  in digital communication,  which maps continuous-amplitude (analog) signals into discrete-amplitude (digital) signals, and hence leads to rounding errors (so called quantization errors) on shared messages. Such quantization errors are usually modeled as additive noises and are non-negligible when the quantizer has a limited number of quantization levels \cite{lee2012digital}. In fact, in deep learning applications where the dimension of optimization variables can scale to hundreds of millions \cite{huang2017densely}, many distributed optimization  algorithms purposely employ a coarse quantizer to reduce the communication overhead \cite{wen2017terngrad,koloskova2019decentralized,cao2021decentralized_submitted}, resulting in large quantization errors. Furthermore, as privacy becomes an increasingly pressing  need  while conventional distributed optimization algorithms are proven to leak information of participating agents \cite{zhang2019admm,zhang2018enabling,melis2019exploiting,wang2022decentralized}, many privacy-aware distributed optimization algorithms opt to inject additive noise in shared messages to ensure differential privacy~\cite{han2016differentially,wang2017differential,zhang2019recycled,he2020differential}. The differential-privacy induced additive noise is persistent throughout the optimization process to ensure a strong privacy protection and will significantly reduce the optimization accuracy of existing distributed optimization algorithms.

In  recent years, adding a decaying factor on the coupling weight has been proven effective in suppressing the influence of persistent information-sharing noise in distributed optimization  \cite{srivastava2011distributed,zhang2018distributed,cao2021decentralized,george2019distributed,doan2021convergence,wang21}. In combination with a decaying stepsize in gradient descent  (to alleviate the effect of information-sharing noise on gradient directions), these approaches can achieve almost sure convergence to an optimal solution. However, these results are only applicable to  static-consensus based distributed optimization algorithms, {which work on symmetric or balanced graphs but cannot be applied to general directed interaction graphs}. In fact, in gradient-tracking based distributed optimization algorithms,  information-sharing noise will accumulate on the estimate of the global gradient, and its variance will grow to infinity when the information-sharing noise is persistent, as will be explained later in Sec. III. Recently,  \cite{pu2020robust} proposed an algorithm which can avoid information-sharing noise from accumulating on the global-gradient estimate when the {inter-agent interaction is} constant. However,  when the {inter-agent interaction is} time-varying, the approach cannot avoid noise accumulation from happening, which precludes the possibility  to incorporate a decaying factor to attenuate the influence of noise. Directly combining conventional gradient-tracking based approaches with a decaying factor
can { reduce the speed of} such noise accumulation but is unable to avoid the {   noise from accumulating on the estimate of the global gradient and the gradient-estimation noise variance from escaping to infinity}. Our recent result exploited heterogeneous decaying factors for the   optimization variable and the gradient-tracking variable, and managed to avoid the accumulated gradient-estimation noise variance from growing to infinity \cite{wang2022tailoring}. However, the gradient-tracking noise   still accumulates with iterations, which significantly affects the accuracy of distributed optimization.
In this paper, we propose to revise the mechanics of the gradient-tracking based approach to tackle information-sharing noise in distributed optimization. More specifically, we propose a new gradient-tracking based architecture which can  avoid noise from accumulation on every agent's estimate of the global gradient. This approach is applicable even when the {inter-agent interaction is} time-varying, which enables the incorporation of  a decaying  factor to attenuate the influence of noise. In fact, by choosing the decaying factor appropriately,   the proposed approach can gradually eliminate the influence of information-sharing noise on all  agents' local gradient-estimates  even when the noise is  persistent, and hence ensures the final optimality of distributed optimization.

The main contributions of this paper are as follows: 1) We propose a new gradient-tracking based distributed optimization architecture that can avoid the accumulation of information-sharing noise in the estimate of the global gradient. Different from \cite{pu2020robust}, which is only applicable when the  { inter-agent interaction is time-invariant}, the new architecture allows  { inter-agent interaction} to be time-varying, which is key  to enable the incorporation of a decaying  factor in the interaction to gradually eliminate the influence of persistent information-sharing noise; 2) By incorporating a decaying  factor in the inter-agent interaction, we arrive at { two new algorithms that are} able to gradually eliminate the influence of information-sharing noise on both consensus and global-gradient estimation, which, to our knowledge, has not been achieved before. {The first algorithm requires each agent to have access to the left eigenvector of the  coupling weight matrix, whereas the second algorithm uses a local eigenvector estimator to avoid requiring such global information}; 3) We prove that even under persistent information-sharing noise, the proposed algorithms can guarantee every agent's almost sure convergence to an optimal solution on general directed graphs. This is in contrast to existing static-consensus based algorithms in \cite{srivastava2011distributed,george2019distributed,doan2021convergence} that are only applicable to  balanced directed  graphs (the sum of each agent's in-neighbor coupling weights equal to the sum of its out-neighbor coupling weights); 4) We prove that the proposed approach can ensure all agents' almost sure convergence to an optimal solution even when the gradient is subject to noise, which is a common problem in machine learning applications;
5) { The proposed convergence analysis has fundamental differences from existing proof techniques for gradient-tracking based algorithms. More specifically, existing  convergence analysis of  gradient-tracking based algorithms  relies  on  formulating  the error dynamics as a linear time-invariant system of inequalities, whose convergence is determined by a constant systems matrix (even under  time-varying coupling graphs \cite{saadatniaki2020decentralized}).  For example, in the existing gradient-tracking based distributed optimization algorithms, this constant systems matrix is denoted as $A$ in \cite{pu2020push,pu2021distributed},    $J$ in \cite{xin2018linear},   $G$ in \cite{xin2019distributed}, or $M$ in \cite{saadatniaki2020decentralized}. Then, existing analysis establishes  exponential (linear) convergence by proving that the spectral radius of  this systems matrix is a constant value strictly less than one. However, under a decaying coupling strength,   the spectral radius of the systems matrix in the conventional formulation will converge to one, which makes it impossible to use the conventional spectral-radius based analysis.  Therefore, to prove convergence of our algorithms, we propose  a new martingale convergence theorem based approach, which is fundamentally different from   conventional proof techniques for gradient-tracking based optimization algorithms.
 Moreover, our algorithms and theoretical derivations only require the objective functions to be convex and Lipschitz continuous in gradients, which is different from many existing results that require the objective functions to be coercive \cite{xu2017convergence} or strongly convex \cite{pu2020robust,xin2019distributed}, or to have bounded gradients \cite{srivastava2011distributed,koloskova2019decentralized,george2019distributed,doan2021convergence,reisizadeh2019exact}.}

The rest of the paper is organized as follows.
Sec.~\ref{se:formulation} formulates the problem
and provides some results for a later use.
Sec.~\ref{se:proposed_algorithm} presents
 a dynamic-consensus based gradient-tracking method that can avoid noise accumulation on the gradient-tracking estimate. Sec. \ref{se:convergence}  establishes the almost sure  convergence of all agents  to a same  optimal solution.  { Sec. \ref{se:estimation_vector} extends the results by incorporating a left-eigenvector estimator into each agent's local update, which ensures a decentralized implementation of the approach even  when  information of the coupling weight matrices is not locally available to individual agents.} Sec.~\ref{se:SGD} extends the results to the case where the gradient is subject to noise and establishes almost sure convergence of all agents to an optimal solution. Sec. \ref{se:simulations} presents numerical comparisons with existing gradient methods  to corroborate the theoretical   results. Finally, Sec.~\ref{se:conclusions} concludes the paper.

{\bf Notations:}  We use $\mathbb{R}^d$ to denote the Euclidean space of
dimension $d$. We write $I_d$ for the identity matrix of dimension $d$,
and ${\bf 1}_d$ for  the $d$-dimensional  column vector will all
entries equal to 1; in both cases  we suppress the dimension when
clear from the context. A vector is viewed as a column vector. For a
vector $x$, $x_i$ denotes its $i$th element.
%We say $x> 0$ (resp. $x\geq 0$) if all elements of
%$x$ are positive (resp. non-negative).
  We use $\langle\cdot,\cdot\rangle$ to denote the inner product of two vectors and
$\|x\|$ for the standard Euclidean norm of a vector $x$.
We write $\|A\|$ for the matrix norm induced by the vector norm $\|\cdot\|$, unless stated otherwise.
We let $A^T$ denote the transpose of a matrix $A$.
We also use other vector/matrix norms defined under a certain transformation determined by a matrix $W$, which will be represented as $\|\cdot\|_W$.
A matrix is column-stochastic when its entries are nonnegative and
elements in every column add up to one. A   matrix $A$ is said
to be row-stochastic when its entries are nonnegative and
elements in every row add up to one. {For two  vectors  $u$ and $v$ with the same dimension, we use $u\leq v$ to represent that every entry of $u$ is no larger than the corresponding entry of $v$.}  Often, we abbreviate {\it almost surely} by {\it a.s}.
\def\as{{\it a.s.\ }}

\section{Problem Formulation and Preliminaries}\label{se:formulation}
We consider  a network of $m$ agents. The agents interact  on a general directed graph. We describe a directed graph using an ordered pair $\mathcal{G}=([m],\mathcal{E})$, where
$[m]=\{1,2,\ldots,m\}$ is the set of nodes (agents) and $\mathcal{E}\subseteq [m]\times [m]$  is the edge set of ordered node pairs describing the interaction among agents.
For a nonnegative weight  matrix $W=\{w_{ij}\}\in\mathbb{R}^{m\times m}$, we define the induced directed graph as $\mathcal{G}_W=([m],\mathcal{E}_W)$, where
the directed edge $(i,j)$ from agent $j$ to agent $i$ exists,
 i.e., $(i,j)\in \mathcal{E}_W$ if and only if $w_{ij}>0$.
For an agent $i\in[m]$,
its in-neighbor set
$\mathbb{N}^{\rm in}_i$ is defined as the collection of agents $j$ such that $w_{ij}>0$; similarly,
the out-neighbor set $\mathbb{N}^{\rm out}_i$ of agent $i$ is the collection of agents $j$ such that $w_{ji}>0$.

By assigning a copy $x_i$ of the decision variable $x$ to each
agent $i$, and then imposing the requirement $x_i = x$ for all
$1 \leq i \leq m$, we can rewrite the optimization problem  (\ref{eq:optimization_formulation1})
  as the following equivalent multi-agent optimization
problem:
\begin{equation}\label{eq:optimization_formulation2}
\min\limits_{x\in\mathbb{R}^{md}}f(x)\triangleq
\frac{1}{m}\sum_{i=1}^m f_i(x_i)\:\: {\rm s.t.}\:\:
x_1=x_2=\cdots=x_m
\end{equation}
where $x_i\in\mathbb{R}^d$ is  agent $i$'s decision variable  and
 the collection
of the agents' variables is
\[x=[x_1^T, x_2^T, \ldots, x_m^T]^T\in\mathbb{R}^{md}.
\]

We make the following standard assumption on the individual objective functions:
\begin{Assumption 1}\label{assumption:f}
Problem (\ref{eq:optimization_formulation1}) has at least one
optimal solution $\theta^{\ast}$. Every $f_i(\cdot)$ is convex and has
Lipschitz continuous gradients, i.e., for some
$L>0$, we have
\[
\|\nabla f_i(u)-\nabla f_i(v)\|\le L \|u-v\|,\quad\forall i \:\:{\rm and}\:\:\forall u,v\in\mathbb{R}^d.
\]
\end{Assumption 1}

\section{The Proposed Approach}\label{se:proposed_algorithm}
In gradient-tracking based algorithms, besides an optimization variable $x_i^k$,
  every agent $i\in[m]$ also maintains and updates a gradient-tracking variable  $y_i^k$ which  estimates   the global gradient (``joint agent"  descent direction).
Both the optimization variable and the gradient-tracking variable  have to be shared with  neighboring agents. The two variables can be shared using two different communication networks, usually called,  $\mathcal{G}_R$ and $\mathcal{G}_C$, which are, respectively, induced by matrices $R\in\R^{m\times m}$ and $C\in\R^{m\times m}$;
that is $(i,j)$ is a directed link in the graph $\mathcal{G}_R$ if and only if $R_{ij}>0$ and, similarly, $(i,j)$ is a directed link in $\mathcal{G}_C$ if and only if $C_{ij}>0$.
We make the following assumption on  $R$ and $C$. (Note that,
given a matrix $A$ with non-negative off-diagonal entries, the induced graph does not depend on the diagonal entries of the matrix. Also,
$\mathcal{G}_{A^T}$ is identical to
$\mathcal{G}_{A}$ with the directions of edges reversed.)

\begin{Assumption 4}\label{Assumption:push_pull topology}
The matrices $R,C\in\R^{m\times m}$ have nonnegative off-diagonal entries
($R_{ij}\geq 0$ and $C_{ij}\geq 0$ for all $i\neq j$). {Their diagonal entries are negative, satisfying
 \begin{equation}\label{eq:diagonal_entries}
  R_{ii}=-\sum_{j\in\mathbb{N}_{R,i}^{\rm in}}R_{ij} ,\quad C_{ii}=-\sum_{j\in\mathbb{N}_{C,i}^{\rm out}}C_{ji}
 \end{equation}
such that $R$ has zero row sums and $C$ has zero column sums. }
 The induced graphs $\mathcal{G}_R$ and
$\mathcal{G}_{C^T}$ satisfy:
\begin{enumerate}
  \item  $\mathcal{G}_R$ is strongly connected, i.e., there is a  path (respecting the directions of edges) from each node to every other node;
  \item  The graph induced by $C^T$, i.e., $\mathcal{G}_{C^T}$,  contains at least one spanning tree.
\end{enumerate}
\end{Assumption 4}

\begin{Remark 1}
The assumption on  $\mathcal{G}_{C^T}$
 is weaker than requiring that the induced graph  $\mathcal{G}_C$ is strongly connected.
\end{Remark 1}

When there is information-sharing noise,   shared messages may be corrupted by  noise. Namely, when agent $i$ shares $x_i^k$ with agent $j$, agent $j$ can only receive a distorted version $x_i^k+\zeta_i^k$ of $x_i^k$, where $\zeta_i^k$ denotes the information-sharing noise. Similarly, when
  agent $i$ shares $y_i^k$ with agent $j$, agent $j$ can only receive a distorted version $y_i^k+\xi_i^k$ of $y_i^k$, where $\xi_i^k$ denotes the information-sharing noise. The  noises $\xi_i^k$ and $\zeta_i^k$ will significantly impact the accuracy of optimization. In fact, as   conventional gradient-tracking algorithms feed  the incremental gradient to the $y$ iterate, the noise on  $y_i^k$ will accumulate and the variance of noise can grow to infinity as iteration proceeds (this will be detailed later).

  To alleviate the influence of information-sharing noise, a decaying factor can be applied to the coupling weight matrix, which has been proven effective in static-consensus based  distributed optimization algorithms \cite{srivastava2011distributed,george2019distributed,doan2021convergence}. However, for gradient-tracking based algorithms, even with a decaying factor on the coupling weight matrices, the noise on $y_i^k$ will still accumulate and increase with time, significantly affecting the accuracy of optimization results.  Recently, \cite{pu2020robust} showed that instead of tracking the global gradient, tracking the cumulative gradient can avoid information noise from accumulating in gradient-tracking based distributed optimization.  However, this approach cannot eliminate the influence of persistent information-sharing noise, and it is subject to steady-state errors. Furthermore, it can only avoid noise accumulation   when the {inter-agent interaction is} time-invariant, precluding the possibility of combining a decaying factor (which will make {inter-agent interaction} time-varying)  to  gradually  attenuate  the influence of noise. In this paper, we propose a new algorithm that can achieve both avoidance of  noise-accumulation   and  incorporation of a decaying factor. By sharing the cumulative-gradient estimate (denoted as an $s$ variable) instead of the direct gradient estimate (i.e., the $y$ variable), we can gradually annihilate the influence of information-sharing noise on the estimate of the global gradient, even when the information-sharing noise is persistent.

\noindent\rule{0.49\textwidth}{0.5pt}
\noindent\textbf{Algorithm 1: Robust gradient-tracking
based distributed optimization}

\vspace{-0.2cm}\noindent\rule{0.49\textwidth}{0.5pt}
\begin{enumerate}[wide, labelwidth=!, labelindent=0pt]
    \item[] Parameters: Stepsize $\lambda^k$ and a decaying factor $\gamma^k$  to suppress information-sharing noise;
    %,  in-bound mixing/pulling weighs $R_{ij}>0$ for all $j\in\mathbb{N}_{R,i}^{\rm in}$, and out-bound pushing weights $C_{l,i}>0$ for all $l\in\mathbb{N}_{C,i}^{\rm out}$, otherwise $R_{ij}=C_{li}=0$;
    \item[] Every agent $i$ maintains two states  $x_i^k$ and
    $s_i^k$, which are initialized randomly with $x_i^0\in\mathbb{R}^d$ and $s_i^0\in\mathbb{R}^d$.
    \item[] {\bf for  $k=1,2,\cdots$ do}
    \begin{enumerate}
        %\item Every agent $j$ injects zero-mean DP-noises $\zeta_j^k$ and $\xi_j^k$  to its states      $x_j^k$ and $y_j^k$, respectively.
        \item Agent $i$ pushes $s_i^k$ to each agent
        $l\in\mathbb{N}_{C,i}^{\rm out}$, which will be received as $s_i^k+\xi_i^k$ due to   information-sharing noise.  And agent $i$ pulls $x_j^k$ from each $j\in\mathbb{N}_{R,i}^{\rm in}$, which will be received as $x_j^k+\zeta_j^k$ due to information-sharing noise. Here the subscript  $R$ or $C$ in neighbor sets indicates the neighbors with respect to  the graphs induced by these matrices.
         \item Agent $i$ chooses $\gamma^k>0$  satisfying
        $1+\gamma^kR_{ii}>0$ and $1+\gamma^kC_{ii}>0$ with $R_{ii}$ and $C_{ii}$ defined in (\ref{eq:diagonal_entries}).
        Then, agent $i$ updates its states as follows:
        \begin{equation}\label{eq:update}
        \begin{aligned}
s_i^{k+1}=&(1+\gamma^kC_{ii})s_i^k+\gamma^k\sum_{j\in \mathbb{N}_{C,i}^{\rm
            in}}C_{ij}(s_j^k+\xi_j^k)\\
            &+ \lambda^k\nabla f_i(x_i^{k}),\\
x_i^{k+1}=&(1+\gamma^k R_{ii})x_i^k+\gamma^k\sum_{j\in \mathbb{N}_{R,i}^{\rm
            in}}R_{ij}(x_j^k+\zeta_j^k)\\
            &-\frac{s_i^{k+1}-s_i^{k}}{u_i},
        \end{aligned}
         \end{equation}
         where $u_i$ denotes the $i$th element of the left eigenvector $u^T$ of $I+\gamma^kR$ associated with eigenvalue 1\footnote{Under Assumption \ref{Assumption:push_pull topology},
  the matrix $I+\gamma^kR$ always has a unique
 positive left eigenvector $u^T$ (associated with eigenvalue 1) satisfying $u^T{\bf 1}=m$ (see details in Lemma \ref{lemma:left_right_eigenvectors}). When $R$ is balanced, $u$ becomes the vector ${\bf 1}$ \cite{horn2012matrix}.}.

                  \item {\bf end}
    \end{enumerate}
\end{enumerate}
\vspace{-0.1cm} \rule{0.49\textwidth}{0.5pt}

{
\begin{Remark 1}
As discussed before, a key difference  between the proposed Algorithm 1 and the existing gradient-tracking based algorithms is that Algorithm 1 introduces a decaying factor $\gamma^k$ to suppress the information-sharing noise. Introducing the decaying factor is reasonable for the following reasons: In the early stages of the iteration, the decaying factor is still far from zero, and hence its attenuation effect on information-sharing is not significant, which  allows the necessary  mixture of information and hence consensus of individual agents' optimization variables; As the iteration proceeds and  individual agents' optimization variables converge to each other (thus diminishing  the need for information-sharing), the decaying factor approaches   zero and hence its attenuation effect on information-sharing noise becomes more severe, which effectively eliminates the influence of information-sharing noise. Of course, to ensure that necessary gradient descent steps and information-mixture operations can be performed, the decaying factor has to decrease slower than $\lambda^k$, which will be specified later in the convergence analysis.
\end{Remark 1}
}

To compare our algorithm with  conventional  gradient-tracking based algorithms, we write the algorithm in   matrix form. Defining
\[
{\bf x}^{k}=\left[\begin{array}{c}(x_1^k)^T\\(x_2^k)^T\\\vdots\\(x_m^k)^T\end{array}\right]\in\mathbb{R}^{m\times d},\:
{\bf s}^{k}=\left[\begin{array}{c}(s_1^k)^T\\(s_2^k)^T\\\vdots\\(s_m^k)^T\end{array}\right]\in\mathbb{R}^{m\times d},
\]
\[
{\bf g}^{k}=\left[\begin{array}{c}(g_1^k)^T\\(g_2^k)^T\\\vdots\\(g_m^k)^T\end{array}\right]\in\mathbb{R}^{m\times d},
\]
with $g_i^k=\nabla f_i(x_i^{k})$ and
\[{\pmb \zeta}_w^{k}=\left[\begin{array}{c}(\zeta_{w1}^k)^T\\(\zeta_{w2}^k)^T\\\vdots\\(\zeta_{wm}^k)^T\end{array}\right]\in\mathbb{R}^{m\times d},\:
{\pmb \xi}_w^{k}=\left[\begin{array}{c}(\xi_{w1}^k)^T\\(\xi_{w2}^k)^T\\\vdots\\(\xi_{wm}^k)^T\end{array}\right]\in\mathbb{R}^{m\times d},\:
\]
   with
\[
 \zeta_{wi}\triangleq \sum_{j\in\mathbb{N}_{R,i}^{\rm in}}R_{ij}\zeta_j^k,\quad
 \xi_{wi}\triangleq \sum_{j\in\mathbb{N}_{C,i}^{\rm in}}C_{ij}\xi_j^k,
\]
we write the dynamics of (\ref{eq:update})  in the following more compact form:
\begin{equation}\label{eq:push-pull-R}
\begin{aligned}
{\bf s}^{k+1}&= C^k  {\bf s}^k +\gamma^k {\pmb\xi}_w^k +\lambda^k{\bf g}^{k}\\
{\bf x}^{k+1}&= R^k  {\bf x}^k+\gamma^k {\pmb \zeta}_w^k- U^{-1} ({\bf s}^{k+1}-{\bf s}^{k})
\end{aligned}
\end{equation}
where
\[R^k=I+\gamma^kR,
\]
\[
C^k=I+\gamma^kC,
\]
and
\[
U={\rm diag}(u_1,\,u_2,\,\cdots,u_m),
\]
with $u_i$ denoting the $i$th element of $R_k$'s  left eigenvector $u$  associated with eigenvalue 1.

It can be seen that in the proposed algorithm, ${\bf s}^k-{\bf s}^{k-1}$ is fed into the optimization variable and acts as the global-gradient estimate. This new approach will avoid the accumulation of information-sharing noise on the global-gradient estimate, which plagues  existing gradient-tracking based approaches. To see this, we use the Push-Pull gradient-tracking algorithm   as an example. In the absence of information-sharing noise, the conventional Push-Pull algorithm takes the following form \cite{pu2020push}:
\begin{equation}\label{eq:push-pull-original}
\begin{aligned}
{\bf x}^{k+1}&= R^k  {\bf x}^k-\lambda^k  {\bf y}^k \\
{\bf y}^{k+1}&= C^k  {\bf y}^k +{\bf g}^{k+1}- {\bf g}^{k}.
\end{aligned}
\end{equation}
 By setting ${\bf y}^{0}={\bf g}^{0}$, one can obtain by induction that
 \[
 {\bf 1}^T {\bf y}^{k}=  {\bf 1}^T {\bf g}^{k},
  \]
  i.e., the agents can track the average gradient $\frac{{\bf 1}^T {\bf g}^{k}}{n}$ by ensuring the consensus of all $y_i^k$ (which leads to $y_i^k=\frac{{\bf 1}^T {\bf y}^{k}}{m}$  for all $i$).

However, when exchanged messages are subject to noises, i.e., exchanged $x_i^k$ and $y_i^k$ are received as $x_i^k+\zeta_i^k$ and $y_i^k+\xi_i^k$, respectively, the update of the conventional Push-Pull becomes (after incorporating a decaying factor $\gamma^k$)
\begin{equation}\label{eq:push-pull-noise}
\begin{aligned}
{\bf x}^{k+1}&= R^k  {\bf x}^k+\gamma^k {\pmb \zeta}_w^k-\lambda^k  {\bf y}^k, \\
{\bf y}^{k+1}&= C^k  {\bf y}^k +\gamma^k {\pmb\xi}_w^k+{\bf g}^{k+1}- {\bf g}^{k},
\end{aligned}
\end{equation}
and one can obtain by induction that
\begin{equation}\label{eq:error_accumulation}
{\bf 1}^T {\bf y}^{k}=  {\bf 1}^T \left({\bf g}^{k}+\sum_{l=0}^{k-1}\gamma^l {\pmb\xi}_w^l\right)
\end{equation} even under ${\bf y}^{0}={\bf g}^{0}$.

 Therefore,
under the conventional Push-Pull algorithm, the
information-sharing noise accumulates with time (even with a decaying factor $\gamma^k$) in the estimate of the global gradient, which significantly
compromises optimization accuracy. (This statement is corroborated by the numerical simulation result  for the conventional Push-Pull algorithm in \cite{pu2020push} in Fig. 1, whose optimization-error variance grows with iterations.) It can be easily verified that other gradient-tracking based distributed optimization algorithms have the same issue of accumulating information-sharing noise.

The proposed algorithm successfully circumvents this problem. In fact, using  the update rule of ${\bf s}^k$ in (\ref{eq:push-pull-R}), one has
\begin{equation}\label{eq:bar_y}
\begin{aligned}
{\bf 1}^T  ({\bf s}^{k+1}-{\bf s}^{k})&=  {\bf 1}^T \left(C^k  {\bf s}^k +\gamma^k {\pmb\xi}_w^k +\lambda^k{\bf g}^{k}-{\bf s}^{k}\right)\\
&={\bf 1}^T  \left( \gamma^kC   {\bf s}^k +\gamma^k {\pmb\xi}_w^k +\lambda^k{\bf g}^{k}\right)\\&= {\bf 1}^T  \left(  \gamma^k {\pmb\xi}_w^k +\lambda^k{\bf g}^{k}\right),
\end{aligned}
\end{equation}
where we used the property ${\bf 1}^TC=0$ from the definition of $C_{ii}$ in (\ref{eq:diagonal_entries}). It is clear that the proposed algorithm avoids information-sharing noise from accumulating  on the gradient estimate.
%It is worth noting that \cite{pu2020robust} also proposed an algorithm to avoid  noise accumulation in the gradient estimate. However, the algorithm in \cite{pu2020robust} is  only applicable when the coupling matrices are constant (time-varying coupling matrices therein will result in noise accumulation in the estimate of the global gradient).
It is worth noting that the proposed algorithm achieves avoidance of   noise-accumulation  even when the inter-agent interaction is time-varying, which enables the incorporation of the decaying factor $\gamma^k$ and further the final elimination of the influence of information-sharing noise on gradient estimate, even when the noises $\zeta_i^k$ and $\xi_i^k$ are persistent. In fact, we can prove that when the decaying factor $\gamma^k$ is chosen appropriately, the proposed algorithm can guarantee that all agents' $x_i^k$ will converge to the same optimal solution almost surely.

\section{Convergence Analysis}\label{se:convergence}
For the convenience of convergence analysis, we first present the following properties for the inter-agent coupling ${R}^k=I+\gamma^kR$ and ${C}^k=I+\gamma^kC$:
\begin{Lemma 1}\label{lemma:left_right_eigenvectors} \cite{horn2012matrix} (or Lemma 1 in \cite{pu2020push})
Under Assumption \ref{Assumption:push_pull topology},
for every $k$, the matrix $I+\gamma^kR$ has a unique
  positive left eigenvector $u^T$ (associated with eigenvalue 1) satisfying $u^T{\bf 1}=m$, and   the matrix $I+\gamma^kC$ has a unique   nonnegative right eigenvector $v$ (associated with eigenvalue 1) satisfying $ {\bf 1}^Tv=m$.
\end{Lemma 1}

{
\begin{Remark 1}\label{re:eigenvector_time_invariant}
It is worth noting that the left eigenvector $u^T$ in Lemma \ref{lemma:left_right_eigenvectors} is time-invariant and independent of $\gamma^k$. In fact, using the definition of left eigenvector, it can be seen that $u^T$ satisfies $u^T(I+\gamma^kR)=u^T$, and thus $u^T(\gamma^kR)=0$ and further $u^TR=0$. Namely, $u^T$ corresponds to the left eigenvector of $R$ associated with eigenvalue $0$. Given that $R$ has zero row-sums according to Assumption \ref{Assumption:push_pull topology}, we know that such a $u^T$ always exists. Similarly, we know that the right eigenvector $v$ of $I+\gamma^kC$ is also time-invariant and independent of $\gamma^k$.
\end{Remark 1}

}

{According to Lemma 3 in~\cite{pu2020push}, we further know that the spectral radius of $\bar{R}^k\triangleq I+\gamma^k R-\frac{ {\bf 1}u^T}{m}$ is equal to $1-\gamma^k|\nu_R|<1$, where $\nu_R$ is an eigenvalue of $R$. Furthermore,   there exists a vector norm $\|{\bf x}\|_R\triangleq\|\tilde{R}{\bf x}\|_2$
(where $\tilde{R}$ is determined by $R$~\cite{pu2020push}) such that $\|\bar{R}^k\|_R<1$
is arbitrarily close to the spectral radius of $\bar{R}^k$, i.e., $1-\gamma^k|\nu_R|<1$. Without loss of generality, we represent this norm as  $\|\bar{R}^k\|_R=1-\gamma^k \rho_R <1$, {where $\rho_R$ is an  arbitrarily close approximation of $|\nu_R|$}. { (Note that for the convergence analysis, we only need the fact that  such an $\tilde{R}$ exists, but do not require knowledge of its explicit expression. For an arbitrarily small difference $\epsilon>0$ between $\|\bar{R}^k\|_R$ and the spectral radius of $\bar{R}^k$,  Lemma 5.6.10 in \cite{horn2012matrix}  provides a constructive way of finding  $\tilde{R}$.  Also see Lemma 5 of the extended version of \cite{pu2020robust}  for more discussions about $\tilde{R}$.)} Similarly,   we have that the spectral radius of $\bar{C}^k\triangleq I+\gamma^k C-\frac{ v{\bf 1}^T}{m}$ is equal to $1-\gamma^k|\nu_C|<1$, where $\nu_C$ is an eigenvalue of $C$. Furthermore,   there exists a vector norm $\|{\bf x}\|_C\triangleq\|\tilde{C}{\bf x}\|_2$   (where $\tilde{C}$ is  determined by $C$~\cite{pu2020push}) such that $\|\bar{C}^k\|_C<1$
is arbitrarily close to the spectral radius of $\bar{C}^k$, i.e., $1-\gamma^k|\nu_C|<1$. Without loss of generality, we represent this norm as  $\|\bar{C}^k\|_C=1-\gamma^k \rho_c <1$, {where $\rho_c$ is an  arbitrarily close approximation of $|\nu_C|$}.}

 For convenience in analysis, we also define  the following { (weighted)} average vectors:
\begin{equation}\label{eq:definiton1}
   \bar{x}^k=\frac{u^T{\bf x}^{k}}{m}, \quad \bar{s}^k=\frac{{\bf 1}^T{\bf s}^{k}}{m},  \quad  \bar{g}^k=\frac{{\bf 1}^T{\bf g}^{k}}{m},
\end{equation}
and
\begin{equation}\label{eq:definition2}
 \bar\zeta_w^k=\frac{{u}^T{\pmb\zeta}_w^k}{m}, \quad  \bar\xi_w^k=\frac{{\bf 1}^T{\pmb\xi}_w^k}{m}.
\end{equation}

To analyze the convergence of the proposed algorithm, we first present a generic convergence result for gradient-tracking based distributed optimization algorithms. To this end, we first  define a matrix norm for  ${\bf x}^k$ following \cite{pu2020push}:
%need to define a $\|\cdot\|_R$ norm based measure of the distance between all local iterates $x_^k, \, x_2^k, \cdots, x_m^k$ and their average $\bar{x}^k$.
\begin{equation}\label{eq:matrix_norm_R}
\|{\bf x}^k\|_R=\left\|\left[\|{\bf x}^k_{(1)}\|_R,\,\|{\bf x}^k_{(2)}\|_R,\cdots,\|{\bf x}^k_{(d)}\|_R \right] \right\|_2
\end{equation}
 where the subscript $2$ denotes the $2-$norm and ${\bf x}^k_{(i)}$ denotes the $i$th column of ${\bf x}^k$ for $1\leq i\leq d$. One can easily see that $\|{\bf x}^k-{\bf 1}\bar{x}^k \|_R$ measures the distance between all $x_i^k$ and their weighted average $\bar{x}^k$.

 Similarly, we  define a matrix norm $\|\cdot\|_C$ for $ {\bf s}^k\triangleq \left[s_1^k,\,s_2^k,\cdots,s_m^k\right]^T\in\mathbb{R}^{m\times d}$:
\begin{equation}\label{eq:matrix_norm_C}
\|{\bf s}^k\|_C=\left\|\left[\|{\bf s}^k_{(1)}\|_C,\,\|{\bf s}^k_{(2)}\|_C,\cdots,\|{\bf s}^k_{(d)}\|_C \right] \right\|_2
\end{equation}
and use $\|{\bf s}^k-v{ \bar{s}}^k \|_C$   to measure  the distance between all $s$ iterates and their   average  $\bar{s}^k$ (weighed by $v$).

We also need the following lemmas about sequences of random vectors:
\begin{Lemma 1}(\cite{wang2022tailoring}, Lemma 4)\label{th-dsystem}
Let  $\{\bv^k\}\subset \mathbb{R}^d$
and $\{\bu^k\}\subset \mathbb{R}^p$ be random nonnegative
vector sequences, and $\{a^k\}$ and $\{b^k\}$ be random nonnegative scalar sequences   such that
\[
\mathbb{E}\left[\bv^{k+1}|\mathcal{F}^k\right]\le (V^k+a^k{\bf 1}{\bf1}^T)\bv^k +b^k{\bf 1} -H^k\bu^k
\]
holds almost surely for all $k\geq 0$, where $\{V^k\}$ and $\{H^k\}$ are random sequences of
nonnegative matrices and
$\mathbb{E}\left[\bv^{k+1}|\mathcal{F}^k \right]$ denotes the conditional expectation given
 $\bv^\ell,\bu^\ell,a^\ell,b^\ell,V^\ell,H^\ell$ for $\ell=0,1,\ldots,k$.
Assume that $\{a^k\}$ and $\{b^k\}$ satisfy
$\sum_{k=0}^\infty a^k<\infty$ and $\sum_{k=0}^\infty b^k<\infty$ almost surely, and
that there exists a (deterministic) vector $\pi>0$ such that
\[\pi^T V^k\le \pi^T,\qquad \pi^TH^k\ge 0,\qquad\forall k\geq 0\]
hold almost surely.
Then, we have
\begin{enumerate}
  \item  $\{\pi^T\bv^k\}$ converges almost surely to some random variable $\pi^T\bv\geq 0$;
  \item  $\{\bv^k\}$ is bounded almost surely;
  \item  $\sum_{ k=0 }^\infty \pi^TH^k\bu^k<\infty$ holds almost surely.
\end{enumerate}
\end{Lemma 1}

%\begin{proof}
%By multiplying the given relation for $\bv^{k+1}$ with $\pi$ and using $\pi^TV^k\le\pi^T$ and the nonnegativity of  $\bv^k$,
%we obtain
%\[
%\mathbb{E}\hspace{-0.07cm}\left[\pi^T\bv^{k+1}|\mathcal{F}^k\right]\hspace{-0.07cm}\le\hspace{-0.07cm} \pi^T\bv^k\hspace{-0.07cm} +a^k(\pi^T {\bf 1})({\bf1}^T\bv^k) +b^k\pi^T{\bf 1}-\pi^TH^k\bu^k
%\]
%Since $\pi>0$, we have $\pi_{\min}=\min_i\{\pi_i\}>0$, which yields
%\[{\bf1}^T\bv^k=\frac{1}{\pi_{\min}} \,\pi_{\min}{\bf1}^T\bv^k\le
%\frac{1}{\pi_{\min}} \,\pi^T\bv^k\]
%where the inequality holds since $\bv^k\ge0$.
%Therefore, one obtains
%\[\mathbb{E}\left[\pi^T\bv^{k+1}|\mathcal{F}^k\right]\le \left(1+a^k\frac{\pi^T {\bf 1}}{\pi_{\min}}\right)\pi^T\bv^k +b^k\pi^T{\bf 1} -\pi^TH^k\bu^k\]
%By our assumption, we have  $\pi^TH^k\bu^k\ge0$ for all $k$ almost surely.
%Thus, the preceding relation implies that
%the conditions of Lemma~\ref{Lemma-polyak} are satisfied
%with $v^k = \pi^T\bv^k$, $\a^k = a^k\pi^T {\bf 1}/\pi_{\min}$ and $\b^k = b^k \pi^T{\bf 1}$.
%So, by Lemma~\ref{Lemma-polyak},
%it follows that the limit $\lim_{k\to\infty}\pi^T\bv^k$ exists almost surely.
%Consequently, $\{\pi^T\bv^k\}$ is bounded almost surely, and
%since $\{\bv^k\}$ is nonnegative  and $\pi>0$, it follows that $\{\bv^k\}$ is also bounded almost surely.
%By Lemma \ref{Lemma-polyak}, we have
%$\sum_{k=0}^\infty \pi^TH^k\bu^k<\infty$  almost surely.
%\end{proof}

\begin{Lemma 1}(\cite{wang2022tailoring}, Lemma 7)\label{Theo:convergence_to_zero}
Let $\{\bv^k\}\subset \mathbb{R}^d$  be a sequence of non-negative random vectors and
 $\{b^k\}$ be a sequence of nonnegative random scalars such that
$\sum_{k=0}^\infty b^k<\infty$
and
\[
\mathbb{E}\left[
\bv^{k+1}|\mathcal{F}^k\right]\le V^k\bv^k +b^k{\bf 1}, \quad \forall k\ge0
\]
hold almost surely,
where $\{V^k\}$ is a sequence of
non-negative matrices  and $\mathcal{F}^k=\{\bv^\ell,b^\ell;0\le \ell\le k\}$.
%, and $\{b^k\}$ is a non-negative scalar sequence such that $\lim_{k\to\infty} b^k=0$.
Assume that there exist a vector $\pi>0$ and a deterministic scalar sequence $\{a^k\}$ satisfying
$a^k\in(0,1)$, $\sum_{k=0}^\infty a^k=\infty$,
and
$\pi^T V^k\le (1-a^k) \pi^T$  for all $k\ge0$.
Then, we have $\lim_{k\to\infty}\bv^k=0$ almost surely.
\end{Lemma 1}
%\begin{proof}
%By multiplying the given relation for $\bv^{k+1}$ with $\pi$ and using
%$\pi^TV^k\le(1-a^k)\pi^T$, we obtain the following relation due to
%the nonnegativity of the vectors $\bv^k$:
%\[
%\mathbb{E}\left[\pi^T\bv^{k+1}|\mathcal{F}^k\right]\le (1-a^k)\pi^T\bv^k +b^k\pi^T{\bf 1}, \quad \forall k\ge0\ \as
%\]
% Since $\sum_{k=0}^\infty a^k= \infty$, and $\sum_{k=0}^\infty b^k<\infty$ {\it a.s.},
%% $\lim_{k\to\infty}b^k/\a^k=0$,
%the conditions of Lemma~\ref{Lemma-polyak_2} are satisfied
%with $v^k = \pi^T\bv^k$, $\a^k=0$, $q^k=a^k$, and $p^k = b^k \pi^T{\bf 1}$.
%By Lemma~\ref{Lemma-polyak_2},
%$\lim_{k\to\infty}\pi^T\bv^k=0$ {\it a.s.} The sequence $\{\bv^k\}$ being nonnegative and $\pi>0$ imply
%that $\lim_{k\to\infty}\bv^k=0$ \as
%\end{proof}

Now we are in a position to present the generic convergence result for gradient-tracking based distributed optimization algorithms:

\begin{Theorem 1}\label{Theorem:general_gradient_tracking}
Assume that the objective function $F(\cdot)$ is continuously
differentiable and that the problem~(\ref{eq:optimization_formulation1}) has an optimal solution $\theta^{\ast}$.
Suppose that a distributed algorithm generates
a sequence  $\{x_i^k\}\subseteq\mathbb{R}^d$ under coupling weight matrix $ R$ and a sequence
$\{s_i^k\}\subseteq\mathbb{R}^d$  under coupling weight matrix $ C$,   such that the following relationship holds
almost surely for some sufficiently large integer $T\ge0$ and for all $k\ge T$:
\begin{equation}\label{eq-fine}
\mathbb{E}\left[\bv^{k+1}|\mathcal{F}^k\right] \le \left(V^k\hspace{-0.07cm} + a^k {\bf 1}{\bf 1}^T\right)\bv^{k}+b^k{\bf
1} - H^k \hspace{-0.07cm}\left[\begin{array}{c}
 \|\nabla F(\bar x^k)\|^2\cr
 \|\bar {g}^k\|^2\end{array}\right]
\end{equation}
 where
 \[
 \mathcal{F}^k=\{x_i^\ell,s_i^\ell;0\le \ell\le k,\ i\in[m]\},
 \]
 and
\[
\begin{aligned}
&\bv^k=\left[\begin{array}{c}\bv^k_1\\\bv^k_2\\\bv^k_3\end{array}\right]\triangleq\left[\begin{array}{c}  F(\bar
x^{k})-F(\theta^{\ast}) \cr  \|{\bf x}^{k}-{\bf 1}\bar x^{k}\|_R^2\cr
\|{\bf s}^{k}- v{ \bar s}^{k}\|_C^2\end{array}
 \right],\\
&
 V^k=
  \left[\begin{array}{ccc} 1 &  \kappa_1\lambda^k &0\\
                            0 & 1-\kappa_2\gamma^k&\kappa_3\gamma^k\\
                            0 & 0&1-\kappa_4\gamma^k
 \end{array}\right],\\
 &H^k=
 \left[\begin{array}{cc} \kappa_5\lambda^k & \kappa_6 \lambda^k-\kappa_7(\lambda^k)^2\\0 &
0\cr 0 & 0
\end{array} \right],
\end{aligned}
\]
with  $\kappa_i>0$ for all  $1\leq i\leq 7$  and $\kappa_2,\kappa_4\in(0,1)$, while the nonnegative scalar
sequences $\{a^k\}$, $\{b^k\}$, and positive sequences  $\{\lambda^k\}$ and $\{\gamma^k\}$ satisfy
$\sum_{k=0}^\infty a^k<\infty$ {\it a.s.}, $\sum_{k=0}^\infty b^k<\infty$ {\it a.s.}, $\sum_{k=0}^\infty \lambda^k=\infty$, $\sum_{k=0}^\infty \gamma^k=\infty$, $\sum_{k=0}^\infty (\gamma^k)^2<\infty$,   $\sum_{k=0}^\infty \frac{(\lambda^k)^2}{\gamma^k}<\infty$, and $\lim_{k\to\infty}\lambda^k/\gamma^k=0$.
Then, we have:
\begin{itemize}
\item[(a)]  $\lim_{k\to\infty} F(\bar x^k)$
exists almost surely and
\[
\lim_{k\to\infty}\|x_i^k - \bar x^k\|= \lim_{k\to\infty}\|s_i^k -
v_i\bar s^k\|=0, \: \forall i,\quad a.s. \]
\item[(b)]
$\displaystyle\lim\inf\limits_{k\rightarrow\infty}\|\nabla F(\bar x^k)\|=0$ holds almost surely.
 Moreover, if the function $F(\cdot)$ has bounded level sets, then
$\{\bar x^k\}$ is bounded and every accumulation point of  $\{\bar
x^k\}$ is an optimal solution almost surely,
and
\[\lim_{k\rightarrow\infty}F(x_i^k)= F(\theta^*),\qquad \forall  i\in[m], \quad \as\]
\end{itemize}

\end{Theorem 1}
\begin{proof}
Since the results of Lemma~\ref{th-dsystem} are asymptotic, they remain valid when the starting index is shifted from $k=0$ to $k=T$, for an arbitrary $T\ge 0$. So the  idea is to show that the conditions in Lemma~\ref{th-dsystem}
are satisfied  for all $k\ge T$, where $T\geq 0$ is large enough.

\noindent(a)  Because $\kappa_i>0$ for all $1\leq i\leq 7$, for $\pi=[\pi_1, \pi_2, \pi_3]^T$ to satisfy  $\pi^T V\leq \pi^T$ and $\pi^TH^k\geq 0$,  we only need to show that the following inequalities are true
\begin{equation}\label{eq:conditons_pi}
\begin{aligned}
 &\kappa_1\lambda^k\pi_1+ (1-\kappa_2\gamma^k) \pi_2 \leq \pi_2,\\
& \kappa_3\gamma^k\pi_2+(1-\kappa_4 \gamma^k) \pi_3\leq \pi_3,
\\
 &\left(\kappa_6 \lambda^k-\kappa_7(\lambda^k)^2\right)\pi_1 \geq 0.
\end{aligned}
\end{equation}

The first inequality is equivalent to $\pi_2\geq \frac{\kappa_1\lambda^k}{\kappa_2\gamma^k}\pi_1$. Given that $\lim_{k\to\infty}\lambda^k/\gamma^k=0$ holds and $\gamma^k$ as well as $\lambda^k$ is positive   according to the assumption, it can easily be seen that for any given $\pi_1>0$, we can always find a $\pi_2>0$ satisfying the relationship when $k$ is larger than some $T \geq 0$.

The second inequality is equivalent to $\pi_3\geq\frac{\kappa_3}{\kappa_4}\pi_2$, which can always be satisfied by setting $\pi_3=\frac{\kappa_3}{\kappa_4}\pi_2$ after fixing $\pi_2$.

The third inequality is equivalent to $\kappa_6  -\kappa_7 \lambda^k >0$, which is always satisfied given that $(\lambda^k)^2$ is summable (and hence $\lambda^k$ tends to zero).

 Thus,   we can always find a vector $\pi$  satisfying all  inequalities in (\ref{eq:conditons_pi}) for $k \ge T$ for some large enough $T\geq 0$, and hence the conditions in Lemma~\ref{th-dsystem} are satisfied.

 By Lemma~\ref{th-dsystem}, it follows that
for the three entries of $\bv^k$, i.e., $\bv_1^k$, $\bv_2^k$, and
$\bv_3^k$, we have that
\be\label{eq-limit-exist} \lim_{k\to\infty}
\pi_1\bv^k_1+\pi_2 \bv_2^k+\pi_3\bv_3^k \ee
%\be\label{eq-limit-exist}
%\lim_{k\to\infty} \left(\pi_1(F(\bar x^k) -
%F^*)+\pi_2 \sum_{i=1}^m\|x_i^k-\bar
%x^k\|^2+\pi_3\sum_{i=1}^m\|y_i^k-\bar y^k\|^2 \right)
%\ee
 exists almost surely, and
\[
\sum_{k=0}^\infty\pi^T H^k\bu^k<\infty
\]
holds almost surely with
\[\bu^k= [\|\nabla
F(\bar x^k)\|^2,\ \|\bar g^k\|^2]^T.
\]
 Since $\pi^TH^k$ has the following form
\[
\begin{aligned}
&\pi^TH^k=\left[\kappa_5\lambda^k\pi_1,\,
(\kappa_6 \lambda^k-\kappa_7(\lambda^k)^2)\pi_1\right]
\end{aligned}
\]
and  $(\lambda^k)^2$ is summable,
 one has \be\label{eq-sumfinite}
\sum_{k=0}^\infty\lambda^k\|\nabla F(\bar x^k)\|^2<\infty, \quad
\sum_{k=0}^\infty\lambda^k\|\bar g^k\|^2<\infty, \quad a. s.\ee
 Hence, it follows that
\be\label{eq-to0}
 \| \nabla F(\bar x^k)\|<\Delta_1,\quad
 \|\bar g^k\|<\Delta_2,\quad a.s.
\ee
for some random scalars $\Delta_1>0$ and $\Delta_2>0$ due to the assumption
$\sum_{k=0}^\infty\lambda^k=\infty$.

Now, we focus on proving that both $\bv_2^k= \|{\bf x}^{k}-{\bf 1}\bar x^{k}\|_R^2$ and $\bv_3^k=
\|{\bf s}^{k}-v{ \bar s}^{k}\|_C^2$ converge to
0 almost surely. The  idea is to show that we can apply Lemma~\ref{Theo:convergence_to_zero}.
By focusing on the second and third elements of $\bv^k$, i.e., $\bv_2^k$ and $\bv_3^k$, from
\eqref{eq-fine} we have
\bas \left[\begin{array}{c} \bv_2^{k+1}\cr
\bv_3^{k+1}\end{array} \right] \le \left(\tilde V^k + a^k {\bf 1}{\bf
1}^T\right) \left[\begin{array}{c} \bv_2^k\cr
\bv_3^k\end{array}\right] +\hat b^k{\bf 1},
\eas
where
\[
\hat b^k = b^k+a^k(F(\bar x^k)-F(\theta^{\ast})),
\]
 \[
 \tilde
V^k=\left[\begin{array}{cc}
 1-\kappa_2\gamma^k & \kappa_3\gamma^k\cr
0 & 1-\kappa_4\gamma^k\end{array}\right],
\]
 which can be rewritten as
\ba\label{eq-finer}
\left[\begin{array}{c} \bv_2^{k+1}\cr
\bv_3^{k+1}\end{array} \right] \le \tilde V^k \left[\begin{array}{c}
\bv_2^k\cr \bv_3^k\end{array}\right] +\tilde b^k{\bf 1} \ea
 with
\[
\begin{aligned}
\tilde b^k=& b^k \\
&+a^k\left( F(\bar x^k) - F(\theta^{\ast})  + \|{\bf x}^{k}-{\bf 1}\bar x^{k}\|_R^2+
\|{\bf s}^{k}-v{  \bar s}^{k}\|_C^2\right).
\end{aligned}
\]

To apply Lemma \ref{Theo:convergence_to_zero}, noting that $\g^k$ is not summable, we show that the { inequality}
 $\tilde \pi^T
\tilde V^k\leq (1-\a\gamma^k) \tilde\pi^T$ has a solution in $\tilde\pi=[\pi_2,\pi_3]$  with
$\pi_2,\pi_3>0$ and $\a\in(0,1)$.

From
\[
\tilde \pi^T
\tilde V^k\leq(1-\a\gamma^k) \tilde \pi^T,
\]
one has
\[
(1-\kappa_2\gamma^k)\pi_2 \leq (1-\alpha\gamma^k)\pi_2
\]
and
\[
 \kappa_3\gamma^k\pi_2+(1-\kappa_4\gamma^k)\pi_3\leq (1-\alpha\gamma^k)\pi_3,
\]
which can be simplified as
$
\alpha\leq \kappa_2$ and $\alpha\leq \kappa_4-\frac{\pi_2}{\pi_3}\kappa_3
$.

Given  $\kappa_2>0$,  $\kappa_3>0$, and $\kappa_4>0$ according to our assumption, we can always find appropriate $\pi_2>0$ and $\pi_3>0$ to make $\alpha\in (0,1)$ hold.

 We next  prove that the condition $\sum_{k=0}^{\infty}\tilde{b}^k<0$ \as
of Lemma \ref{Theo:convergence_to_zero} is
also   satisfied. Indeed, the condition can be met because: (1) $b^k$ and $a^k$ are summable according to the assumption of the theorem; and (2)   $F(\bar x^k) -
F(\theta^{\ast})$, $\|{\bf x}^{k}-{\bf 1} \barx^{k}\|_R^2$, $\|{\bf s}^{k}-v{ \bar s}^{k}\|_C^2$ are all bounded almost surely
due to  the existence of the limit in (\ref{eq-limit-exist}).
Thus, all the conditions of
Lemma~\ref{Theo:convergence_to_zero} are satisfied, and thus
it follows that  $\lim_{k\to\infty}\|x_i^k - \bar
x^k\|=0$ and $\lim_{k\to\infty}\|s_i^k - v_i\bar s^k\|=0$ hold almost surely.
 Moreover, in view of the existence of the limit in (\ref{eq-limit-exist}) and the facts that $\pi_1>0$ and $v_1^k=F(\bar x^k)-F(\theta^*)$, it follows that $\lim_{k\to\infty} F(\bar x^k)$
exists almost surely.

\noindent(b) Since $\sum_{k=0}^\infty\lambda^k\|\nabla F(\bar
x^k)\|^2<\infty$ holds almost surely (see~\eqref{eq-sumfinite}), from $\sum_{k=0}^\infty
\lambda^k =\infty$, it follows that we have
$\displaystyle\lim\inf_{k\to\infty}\|\nabla F(\bar x^k)\|=0$ almost surely.

 Now, if the function $F(\cdot)$ has bounded level sets, then the sequence
$\{\bar x^k\}$ is   bounded almost surely since the limit $\lim_{k\to\infty} F(\bar x^k)$
exists almost surely (as shown in part (a)). Thus, $\{\bar x^k\}$
has accumulation points almost surely.
Let $\{\bar x^{k_i}\}$ be a sub-sequence such that
$\lim_{i\to\infty}\|\nabla F(\bar x^{k_i})\|=0$ holds almost surely.
Without loss of
generality, we may assume that $\{\bar x^{k_i}\}$ is almost surely convergent, for
otherwise we would choose a sub-sequence of $\{\bar x^{k_i}\}$. Let
$\lim_{i\to\infty}\bar x^{k_i}=\hat x$. Then, by the continuity of the
gradient $\nabla F(\cdot)$, it follows that $\nabla F(\hat x)=0$,
implying that $\hat x$ is an optimal point. Since $F(\cdot)$ is continuous,
 we have  $\lim_{i\to\infty} F(\bar x^{k_i}) = F(\hat x)=F(\theta^*)$.
By part (a),  $\lim_{k\to\infty} F(\bar x^k)$
exists almost surely, and thus we must have $\lim_{k\to\infty} F(\bar x^k)=F(\theta^*)$ almost surely.

Finally,  by part (a), we  have $\lim_{k\to\infty}\|x_i^k-\bar
x^k\|^2=0$ almost surely for every $i$. Thus,  each  $\{x_i^k\}$
has the same accumulation points as the sequence $\{\bar x^k\}$ almost surely, implying by
the continuity of the function $F(\cdot)$ that
 $\lim_{k\to\infty}F(x^k_i)=F(\theta^*)$ holds almost surely for all $i$.
\end{proof}

\begin{Remark 1}
In Theorem~\ref{Theorem:general_gradient_tracking}(b),  the bounded level set condition
can be replaced with
any other condition ensuring that the sequence $\{\bar x^k\}$ is bounded almost surely.
\end{Remark 1}

 Theorem~\ref{Theorem:general_gradient_tracking} is critical for establishing
convergence properties of the gradient tracking-based distributed algorithm together
with suitable conditions on the information-sharing noise. We make the following assumption on the noise:
\begin{Assumption 1}\label{assumption:dp-noises-intrack}
For every $i\in[m]$, the noise sequences $\{\zeta_i^k\}$ and
$\{\xi_i^k\}$ are zero-mean independent random variables, and independent of $\{x_i^0;i\in[m]\}$.
Also, for every $k$, the noise collection $\{\zeta_j^k,\xi_j^k; j\in[m]\}$ is independent. The noise variances
$(\sigma_{\zeta,i}^k)^2=\mathbb{E}\left[\|\zeta_i^k\|^2\right]$ and $(\sigma_{\xi,i}^k)^2=\mathbb{E}\left[\|\xi_i^k\|^2\right]$
and the decaying factor  $\g^k$   are
such that
\begin{equation}\label{eq:condition_assumption5}
\begin{aligned}
\hspace{-0.2cm}\sum_{k=0}^\infty(\g^k)^2&\max_{i\in[m]}(\sigma_{\zeta,i}^k)^2<\infty,  \sum_{k=0}^\infty (\gamma^k)^2\, \max_{j\in[m]}(\sigma_{\xi,j}^k)^2 <\infty.
\end{aligned}
\end{equation}
The initial random vectors satisfy
$\mathbb{E}\left[\|x_i^0\|^2\right]<\infty$,  $\forall i\in[m]$.
\end{Assumption 1}
\begin{Remark 1}
The condition (\ref{eq:condition_assumption5}) is satisfied, for example, when sequences
$\{(\gamma^k)^2\}$ and $\{(\lambda^k)^2\}$ are  summable, and sequences $\{\sigma_{\zeta,i}^k\}$ and
$\{\sigma_{\xi,i}^k\}$ are bounded for every $i\in[m]$.
\end{Remark 1}
\begin{Theorem 1}\label{th:deterministic}
 Let  Assumption 1, Assumption~\ref{Assumption:push_pull topology}, and
Assumption~\ref{assumption:dp-noises-intrack} hold.
If  $\{\gamma^k\}$    and $\{\lambda^k\}$ satisfy
$
 \sum_{k=0}^\infty \gamma^k=\infty$,  $\sum_{k=0}^\infty (\gamma^k)^2<\infty$, $\sum_{k=0}^\infty \lambda^k=\infty$, $\sum_{k=0}^\infty \frac{(\lambda^k)^2}{\gamma^k}<\infty$, and $\lim_{k\to\infty}\lambda^k/\gamma^k=0$,
then the results of Theorem~\ref{Theorem:general_gradient_tracking}
hold for Algorithm 1.
\end{Theorem 1}
\begin{proof}

The  goal is to establish the relationship in~(\ref{eq-fine}), with the $\sigma$-field
$\mathcal{F}^k =\{x_i^\ell,s_i^\ell; 0\le \ell\le k,\ i\in[m]\}$.
To this end, we divide the derivations into four
steps: in Step I, Step II, and Step III, we establish  relations
for $\|{\bf s}^k-v{ \bar{s}}^k \|_C$, $\|{\bf x}^k-{\bf 1}{ \bar{x}}^k \|_R$, and $\mathbb{E}\left[ F(\bar x^k)-F(\theta^{\ast})|\mathcal{F}^k\right]$ for the iterates generated by the proposed algorithm, respectively. In Step IV, we use them to
show that~(\ref{eq-fine}) of Theorem~\ref{Theorem:general_gradient_tracking} holds.

Step I: Relationship for $\|{\bf s}^k-v{ \bar{s}}^k \|_C$.

From (\ref{eq:push-pull-R}), we have
\begin{equation}\label{eq:bar_s}
\begin{aligned}
{\bar s}^{k+1}&=\frac{{\bf 1}^T{\bf s}^{k+1}}{m}\\
&=\frac{{\bf 1}^T}{m}\left(C^k  {\bf s}^k +\gamma^k {\pmb\xi}_w^k +\lambda^k{\bf g}^{k}\right)\\
&= {\bar s}^{k}+\gamma^k\bar{\xi}_w^k+\lambda^{k}\bar{g}^{k},
\end{aligned}
\end{equation}
which, in combination with the relationship $\left(C^k -\frac{v {\bf 1}^T}{m}\right)v=0$,   leads to
\[
{\bf s}^{k+1}-v{ \bar{s}}^{k+1} =\bar{C}^k ( {\bf s}^k -v{ \bar{s}}^{k}) +\gamma^k \Pi_v{\pmb\xi}_w^k+\lambda^{k}\Pi_v {\bf g}^{k},
\]
where we  used the relationship $\bar{C}^k=C^k -\frac{v {\bf 1}^T}{m}$ and defined $\Pi_v=I-\frac{v{\bf 1}^T}{m}$ for the sake of notational simplicity.

The preceding relationship further leads to
\begin{equation}
\begin{aligned}
&\left\|{\bf s}^{k+1}-v{ \bar{s}}^{k+1}\right\|^2_C \\
&=\left\|\bar{C}^k ( {\bf s}^k -v{ \bar{s}}^{k})+\lambda^{k}\Pi_v {\bf g}^{k}\right\|_C^2+\left\|\gamma^k \Pi_v{\pmb\xi}_w^k\right\|_C^2\\
&\quad +2\left\langle \bar{C}^k ( {\bf s}^k -v{ \bar{s}}^{k})+\lambda^{k}\Pi_v {\bf g}^{k}, \gamma^k \Pi_v{\pmb\xi}_w^k\right \rangle_C\\
&\leq \left(\|\bar{C}^k\|_C \| {\bf s}^k -v{ \bar{s}}^{k}\|_C+\lambda^{k}\|\Pi_v\|_C\| {\bf g}^{k}\|_C\right)^2\\
&\quad +\left\|\gamma^k \Pi_v{\pmb\xi}_w^k\right\|_C^2\\
&\quad+2\left\langle \bar{C}^k ( {\bf s}^k -v{ \bar{s}}^{k})+\lambda^{k}\Pi_v {\bf g}^{k}, \gamma^k \Pi_v{\pmb\xi}_w^k\right \rangle_C,
\end{aligned}
\end{equation}
 where $\la\cdot\ra_C$ denotes the inner product induced\footnote{Since one can verify that $\|{\bf s}^k\|_C=\|\tilde{C}{\bf s}^k\|_2$ where   $\tilde{C}$ is discussed in the paragraph after Remark \ref{re:eigenvector_time_invariant}, we have  the norm $\|\cdot\|_C$ satisfying the Parallelogram law, implying that it has an associated inner product $\langle\cdot,\cdot\rangle_C$.} by the norm $\|\cdot\|_C$.

We further bound the first term on the right hand side of the preceding inequality using the property $\|\bar{C}^k\|_C=1-\gamma^k\rho_c$ and the inequality $(a+b)^2\le (1+\epsilon) a^2 + (1+\epsilon^{-1})b^2$ valid for any scalars $a$, $b$, and $\epsilon>0$ (by setting $\epsilon= \frac{1}{1-\gamma^k\rho_c}-1$ and hence $1-\epsilon^{-1}=\frac{1}{\gamma^k\rho_c}$):
\[
\begin{aligned}
&\left\|{\bf s}^{k+1}-v{ \bar{s}}^{k+1}\right\|^2_C\\
&\leq  \left((1-\gamma^k\rho_c)\| {\bf s}^k -v{ \bar{s}}^{k}\|_C+\lambda^{k}\|\Pi_v\|_C \|{\bf g}^{k}\|_C\right)^2\\
&\quad +\left\|\gamma^k \Pi_v{\pmb\xi}_w^k\right\|_C^2\\
&\quad +2\left\langle \bar{C}^k ( {\bf s}^k -v{ \bar{s}}^{k})+\lambda^{k}\Pi_v {\bf g}^{k}, \gamma^k \Pi_v{\pmb\xi}_w^k\right \rangle_C\\
&\leq    (1-\gamma^k\rho_c)\| {\bf s}^k -v{ \bar{s}}^{k}\|^2_C\\
&\quad +\frac{(\lambda^{k})^2}{\gamma^k\rho_c}\|\Pi_v\|^2_C \|{\bf g}^{k}\|_C^2+\left\|\gamma^k \Pi_v{\pmb\xi}_w^k\right\|_C^2\\
&\quad +2\left\langle \bar{C}^k ( {\bf s}^k -v{ \bar{s}}^{k})+\lambda^{k}\Pi_v {\bf g}^{k}, \gamma^k \Pi_v{\pmb\xi}_w^k\right \rangle_C.
\end{aligned}
\]

Taking the expectation (conditioned on $\mathcal{F}^k$) on both sides yields
\begin{equation}\label{eq:s_1}
\begin{aligned}
& \mathbb{E}\left[ \left\|{\bf s}^{k+1}-v{ \bar{s}}^{k+1}\right\|^2_C|\mathcal{F}^k\right]\leq (1-\gamma^k\rho_c)\left\| {\bf s}^k -v{ \bar{s}}^{k} \right\|_C^2 \\
  &\quad + \frac{(\lambda^{k})^2}{\gamma^k\rho_c}\|\Pi_v\|^2_C \|{\bf g}^{k}\|_C^2+ (\gamma^k)^2\|\Pi_v\|_C^2   \mathbb{E}\left[ \| {\pmb\xi}_w^k\|_C^2\right]\\
 &\leq (1-\gamma^k\rho_c)\left\| {\bf s}^k -v{ \bar{s}}^{k} \right\|_C^2 +\frac{(\lambda^{k})^2\delta^2_{C,2}}{\gamma^k\rho_c}\|\Pi_v\|^2_C \|{\bf g}^{k}\|_2^2 \\
 &\quad+ (\gamma^k)^2\delta^2_{C,2}\|\Pi_v\|_C^2   \mathbb{E}\left[ \| {\pmb\xi}_w^k\|_2^2\right]\\
 &=(1-\gamma^k\rho_c)\left\| {\bf s}^k -v{ \bar{s}}^{k} \right\|_C^2 +\frac{(\lambda^{k})^2\delta^2_{C,2}}{\gamma^k\rho_c}\|\Pi_v\|^2_C \|{\bf g}^{k}\|_2^2 \\
  &\quad+ (\gamma^k)^2\delta^2_{C,2}\|\Pi_v\|_C^2   \sum_{i,j} (C_{ij}\sigma_{\xi,j}^k)^2,
\end{aligned}
\end{equation}
where $\delta_{C,2}$ is a constant such that $\|x\|_C\leq \delta_{C,2}\|x\|_2$ for all $x$. (In finite-dimensional vector spaces, all norms are equivalent up to a proportionality constant,  represented by $\delta_{C,2}$ here.) { Note that the   inner-product
term  in the preceding step disappears because the means of all $\xi_i^k$ are zero  according to Assumption \ref{assumption:dp-noises-intrack}, and hence their linear combination ${\pmb\xi}_w^k$ also has zero mean.}

Next we proceed to bound the term $\|{\bf g}^k\|_2$ on the right hand side of the preceding  inequality.

Because every $f_i(\cdot)$ is convex with Lipschitz continuous gradient $L$ according to Assumption \ref{assumption:f}, we always have the following relation (see Theorem 2.1.5 in \cite{Nesterov-book}):
\[
f_i(v)+ \langle \nabla f_i(v), u-v\rangle+\frac{\|\nabla f_i(v)-\nabla f_i(u)\|^2}{2 L}\leq f_i(u)
\]
for any $u,v\in\mathbb{R}^d$.

Letting $v=\theta^{\ast}$ and $u=\bar{x}^k$ in the preceding relation, we obtain for all $i$
\[
f_i(\theta^\ast)+\left\langle\nabla f_i(\theta^\ast),\bar{x}^k-\theta^\ast\right\rangle+\frac{\|\nabla f_i(\theta^\ast)-\nabla f_i(\bar{x}^k) \|^2 }{2L}\leq f_i(\bar{x}^k)
\]
and further
\[
\begin{aligned}
&F(\theta^\ast)+\left\langle\nabla F(\theta^\ast),\bar{x}^k-\theta^\ast\right\rangle+\frac{\sum_{i=1}^{m}\|\nabla f_i(\theta^\ast)-\nabla f_i(\bar{x}^k)\|^2}{2mL}\\
&\leq F(\bar{x}^k).
\end{aligned}
\]
Recalling $\nabla F(\theta^\ast)=0$, we have
\[
 \sum_{i=1}^{m}\|\nabla f_i(\theta^\ast)-\nabla f_i(\bar{x}^k)\|^2\leq 2mL( F(\bar{x}^k)-F(\theta^\ast))
\]
and further
\begin{equation}\label{eq:nabla_f_i}
\begin{aligned}
&\sum_{i=1}^{m}\|\nabla f_i(\bar{x}^k)\|^2\\
&\leq   2\sum_{i=1}^{m}\left(\|\nabla f_i(\theta^\ast)-\nabla f_i(\bar{x}^k)\|^2 + \|\nabla f_i(\theta^\ast)\|^2 \right)\\
&\leq           4mL( F(\bar{x}^k)-F(\theta^\ast))+ 2\sum_{i=1}^{m} \|\nabla f_i(\theta^\ast)\|^2.
\end{aligned}
\end{equation}
Therefore, we have
\begin{equation}\label{eq:g_norm}
\begin{aligned}
 \|{\bf g}^k\|^2&=\sum_{i=1}^{m}\|g_i^k\|^2 \\
 &\leq  2\sum_{i=1}^{m}\left(\|g_i^k-\nabla f_i(\bar{x}^k)\|^2 + \|\nabla f_i(\bar{x}^k)\|^2 \right)\\
 &\leq 2L^2\sum_{i=1}^{m} \|x_i^k- \bar{x}^k \|^2  + 8mL( F(\bar{x}^k)-F(\theta^\ast))\\
 &\qquad+ 4\sum_{i=1}^{m} \|\nabla f_i(\theta^\ast)\|^2\\
 &= 2L^2 \|{\bf x}^k- {\bf 1}\bar{x}^k \|_2^2  + 8mL( F(\bar{x}^k)-F(\theta^\ast))\\
 &\qquad+ 4\sum_{i=1}^{m} \|\nabla f_i(\theta^\ast)\|^2.
\end{aligned}
\end{equation}

Plugging (\ref{eq:g_norm}) into (\ref{eq:s_1}) yields

\begin{equation}\label{eq:s_final}
\begin{aligned}
& \mathbb{E}\left[ \left\|{\bf s}^{k+1}-v{ \bar{s}}^{k+1}\right\|^2_C|\mathcal{F}^k\right] \leq(1-\gamma^k\rho_c)\left\| {\bf s}^k -v{ \bar{s}}^{k} \right\|_C^2 \\
&\qquad+\frac{2L^2(\lambda^{k})^2\delta^2_{C,2}\|\Pi_v\|^2_C}{\gamma^k\rho_c}  \|{\bf x}^k- {\bf 1}\bar{x}^k \|_2^2\\
&\qquad+\frac{8mL(\lambda^{k})^2\delta^2_{C,2}\|\Pi_v\|^2_C}{\gamma^k\rho_c} ( F(\bar{x}^k)-F(\theta^\ast))   \\
&\qquad+\frac{4(\lambda^{k})^2\delta^2_{C,2}\|\Pi_v\|^2_C}{\gamma^k\rho_c} \sum_{i=1}^{m} \|\nabla f_i(\theta^\ast)\|^2\\
  &\qquad+ (\gamma^k)^2\delta^2_{C,2}\|\Pi_v\|_C^2   \sum_{i,j} (C_{ij}\sigma_{\xi,j}^k)^2.
\end{aligned}
\end{equation}

Step II: Relationship for $\|{\bf x}^k-{\bf 1}{ \bar{x}}^k \|_R$.

From (\ref{eq:push-pull-R}), we obtain
\begin{equation}\label{eq:bar_x}
\begin{aligned}
\bar{x}^{k+1}&= \frac{u^T}{m}{\bf x}^{k+1}= \frac{u^T}{m}( R^k  {\bf x}^k+\gamma^k {\pmb \zeta}_w^k- U^{-1} ({\bf s}^{k+1}-{\bf s}^{k}))\\
&=\bar{x}^k+\gamma^k\bar\zeta_w^k-\frac{{\bf 1}^T}{m} ({\bf s}^{k+1}-{\bf s}^{k})\\
&=\bar{x}^k+\gamma^k\bar\zeta_w^k-({\bar s}^{k+1}-{\bar s}^{k})\\
&=\bar{x}^k+\gamma^k\bar\zeta_w^k-\gamma^k\bar{\xi}_w^k-\lambda^{k}\bar{g}^{k},
\end{aligned}
\end{equation}
where we used  $u^TU^{-1}={\bf 1}^T$ in the second equality and  (\ref{eq:bar_s}) in the last equality.

Combining (\ref{eq:push-pull-R}) and (\ref{eq:bar_x})  leads to
\begin{equation}\label{eq:x_bar_x}
\begin{aligned}
{\bf x}^{k+1}-{\bf 1}\bar{x}^{k+1}=&\bar{R}^k ({\bf x}^{k}-{\bf 1}\bar{x}^{k})-\Pi_U({\bf s}^{k+1}-{\bf s}^k)+\gamma^k\Pi_u{\pmb \zeta}_w^k,
\end{aligned}
\end{equation}
where we used the relationship $\bar{R}^k{\bf 1}\bar{x}^{k}=0$ and $\bar{R}^k=R^k-\frac{{\bf 1}u^T}{m}$, and defined $\Pi_u=I-\frac{{\bf 1}u^T}{m}$, $\Pi_U=U^{-1}-\frac{{\bf 1}{\bf 1}^T}{m}$ for the sake of notational simplicity.

From the first relationship in (\ref{eq:push-pull-R}), we can obtain
\begin{equation}\label{eq:s_s}
\begin{aligned}
{\bf s}^{k+1}-{\bf s}^k&=C^k  {\bf s}^k +\gamma^k {\pmb\xi}_w^k +\lambda^k{\bf g}^{k}-{\bf s}^k\\
&=\gamma^kC{\bf s}^k+\gamma^k {\pmb\xi}_w^k +\lambda^k{\bf g}^{k}\\
&=\gamma^kC({\bf s}^k-v\bar{s}^k)+\gamma^k {\pmb\xi}_w^k +\lambda^k{\bf g}^{k},
\end{aligned}
\end{equation}
where we used ${C}^k=I+\gamma^kC$ in the second equality and $Cv=0$ in the last equality.

Combining (\ref{eq:x_bar_x}) and (\ref{eq:s_s}) yields
\begin{equation}\label{eq:x_bar_x2}
\begin{aligned}
\begin{aligned}
{\bf x}^{k+1}-{\bf 1}\bar{x}^{k+1}=&\bar{R}^k ({\bf x}^{k}-{\bf 1}\bar{x}^{k})-\gamma^k\Pi_U C({\bf s}^k-v\bar{s}^k)\\
&  -\gamma^k\Pi_U {\pmb\xi}_w^k  -\lambda^k\Pi_U{\bf g}^{k}  +\gamma^k\Pi_u{\pmb \zeta}_w^k.
\end{aligned}
\end{aligned}
\end{equation}

Taking the norm $\|\cdot\|_R$ on both sides of the preceding relationship yields
\begin{equation}\label{eq:norm_x}
\begin{aligned}
&\|{\bf x}^{k+1}-{\bf 1}\bar{x}^{k+1}\|_R^2\\
&=\|\bar{R}^k ({\bf x}^{k}-{\bf 1}\bar{x}^{k})-\gamma^k\Pi_U C({\bf s}^k-v\bar{s}^k)-\lambda^k\Pi_U{\bf g}^{k}\|_R^2 \\
&\quad+ (\gamma^k)^2\|\Pi_u{\pmb \zeta}_w^k-\Pi_U{\pmb \xi}_w^{k} \|_R^2 \\
&\quad +2\left\langle \bar{R}^k ({\bf x}^{k}-{\bf 1}\bar{x}^{k})-\gamma^k\Pi_U C({\bf s}^k-v\bar{s}^k)-\lambda^k\Pi_U{\bf g}^{k},\right.\\
&\qquad \left.\gamma^k\Pi_u{\pmb \zeta}_w^k-\gamma^k\Pi_U{\pmb \xi}_w^{k}\right\rangle_R\\
&\leq \left(\|\bar{R}^k\|_R \|{\bf x}^{k}-{\bf 1}\bar{x}^{k}\|_R+\gamma^k\|\Pi_U C\|_R\|{\bf s}^k-v\bar{s}^k\|_R\right. \\
&\quad \left.+\lambda^k\|\Pi_U\|_R\|{\bf g}^{k}\|_R\right)^2+ (\gamma^k)^2\|\Pi_u{\pmb \zeta}_w^k-\Pi_U{\pmb \xi}_w^{k} \|_R^2 \\
&\quad +2\left\langle \bar{R}^k ({\bf x}^{k}-{\bf 1}\bar{x}^{k})-\gamma^k\Pi_U C({\bf s}^k-v\bar{s}^k)-\lambda^k\Pi_U{\bf g}^{k},\right.\\
&\qquad \left.\gamma^k\Pi_u{\pmb \zeta}_w^k-\gamma^k\Pi_U{\pmb \xi}_w^{k}\right\rangle_R,
\end{aligned}
\end{equation}
 where $\la\cdot\ra_R$ denotes the inner product induced\footnote{Since one can verify that $\|{\bf x}^k\|_R=\|\tilde{R}{\bf x}^k\|_2$ where $\tilde{R}$ is discussed in the paragraph after Remark \ref{re:eigenvector_time_invariant}, we have  the norm $\|\cdot\|_R$ satisfying the Parallelogram law, implying that it has an associated inner product $\langle\cdot,\cdot\rangle_R$.} by the norm $\|\cdot\|_R$.

Using the relationship $\|\bar{R}^k\|_R=1-\gamma^k\rho_R$ and the inequality $(a+b)^2\le (1+\epsilon) a^2 + (1+\epsilon^{-1})b^2$ valid for any scalars $a$, $b$, and $\epsilon>0$ (by setting $\epsilon= \frac{1}{1-\gamma^k\rho_R}-1$ and hence $1-\epsilon^{-1}=\frac{1}{\gamma^k\rho_R}$), we can arrive at
\begin{equation}\label{eq:Algorithm1_x_error}
\begin{aligned}
&\|{\bf x}^{k+1}-{\bf 1}\bar{x}^{k+1}\|_R^2\leq (1-\gamma^k\rho_R) \|{\bf x}^{k}-{\bf 1}\bar{x}^{k}\|_R^2\\
&\quad+\frac{2\gamma^k\|\Pi_UC\|_R^2}{\rho_R}\|{\bf s}^k-v\bar{s}^k\|_R^2 +\frac{2(\lambda^k)^2\|\Pi_U\|^2_R}{\gamma^k\rho_R}\|{\bf g}^{k}\|^2_R\\
&\quad+(\gamma^k)^2\|\Pi_u{\pmb \zeta}_w^k-\Pi_U{\pmb \xi}_w^{k} \|_R^2 \\
&\quad +2\left\langle \bar{R}^k ({\bf x}^{k}-{\bf 1}\bar{x}^{k})-\gamma^k\Pi_U C({\bf s}^k-v\bar{s}^k)-\lambda^k\Pi_U{\bf g}^{k},\right.\\
&\qquad \left.\gamma^k\Pi_u{\pmb \zeta}_w^k-\gamma^k\Pi_U{\pmb \xi}_w^{k}\right\rangle_R.
\end{aligned}
\end{equation}
Taking the expectation (conditioned on $\mathcal{F}^k$) on both sides yields
\begin{equation}\label{eq:x-bar_x_2}
 \begin{aligned}
& \mathbb{E}\left[ \|{\bf x}^{k+1}-{\bf 1}\bar{x}^{k+1}\|_R^2 |\mathcal{F}^k\right]\leq (1-\gamma^k\rho_R) \|{\bf x}^{k}-{\bf 1}\bar{x}^{k}\|_R^2\\
&\quad+\frac{2\gamma^k\|\Pi_UC\|_R^2}{\rho_R}\|{\bf s}^k-v\bar{s}^k\|_R^2 +\frac{2(\lambda^k)^2\|\Pi_U\|^2_R}{\gamma^k\rho_R}\|{\bf g}^{k}\|^2_R\\
&\quad+2(\gamma^k)^2\|\Pi_u\|^2\mathbb{E}\left[\|{\pmb \zeta}_w^k\|_R^2\right]+2(\gamma^k)^2\|\Pi_U\|_R^2\mathbb{E}\left[\|{\pmb \xi}_w^{k} \|_R^2\right]\\
&\leq (1-\gamma^k\rho_R) \|{\bf x}^{k}-{\bf 1}\bar{x}^{k}\|_R^2+\frac{2\gamma^k\|\Pi_UC\|_R^2}{\rho_R} \|{\bf s}^k-v\bar{s}^k\|_R^2 \\
&\quad+
 \frac{2(\lambda^k)^2\|\Pi_U\|_R^2\delta^2_{R,2}}{\gamma^k\rho_R} \|{\bf g}^k\|_2^2 \\
 &\quad+ 2(\gamma^k)^2\|\Pi_u\|_R^2   \delta^2_{R,2} \sum_{i,j} (R_{ij}\sigma_{\zeta,j}^{k})^2   \\
 &\quad+2(\gamma^k)^2\|\Pi_U\|_R^2 \delta^2_{R,2} \sum_{i,j} (C_{ij}\sigma_{\xi,j}^{k})^2,
 \end{aligned}
\end{equation}
where $\delta_{R,2}$ is a constant such that $\|x\|_R\leq \delta_{R,2}\|x\|_2$ for all $x$. (As mentioned earlier, in finite-dimensional vector spaces, all norms are equivalent up to a proportionality constant,  represented here by $\delta_{R,2}$.)

Plugging (\ref{eq:g_norm}) into (\ref{eq:x-bar_x_2}) yields

\begin{equation}\label{eq:x-bar_x_final}
 \begin{aligned}
 \begin{aligned}
& \mathbb{E}\left[ \|{\bf x}^{k+1}-{\bf 1}\bar{x}^{k+1}\|_R^2 |\mathcal{F}^k\right]\\
&\leq \left(1-\gamma^k\rho_R+\frac{4(\lambda^k)^2 L^2\|\Pi_U\|_R^2\delta^2_{R,2}}{\gamma^k\rho_R}\right) \|{\bf x}^{k}-{\bf 1}\bar{x}^{k}\|_R^2\\
&\quad+\frac{2\gamma^k\|\Pi_UC\|_R^2}{\rho_R} \|{\bf s}^k-v\bar{s}^k\|_R^2 \\
 &\quad+
 \frac{16mL(\lambda^k)^2\|\Pi_U\|_R^2\delta^2_{R,2}}{\gamma^k\rho_R} (     F(\bar{x}^k)-F(\theta^\ast)) \\
 &\quad+
 \frac{8(\lambda^k)^2\|\Pi_U\|_R^2\delta^2_{R,2}}{\gamma^k\rho_R}   \sum_{i=1}^{m} \|\nabla f_i(\theta^\ast)\|^2    \\
 &\quad+ 2(\gamma^k)^2\|\Pi_u\|_R^2   \delta^2_{R,2} \sum_{i,j} (R_{ij}\sigma_{\zeta,j}^{k})^2   \\
 &\quad+2(\gamma^k)^2\|\Pi_U\|_R^2 \delta^2_{R,2} \sum_{i,j} (C_{ij}\sigma_{\xi,j}^{k})^2.
 \end{aligned}
 \end{aligned}
\end{equation}

Step III: Relationship for $F(\bar{x}^{k})-F(\theta^{\ast})$.

Because $F(\cdot)$ is convex with Lipschitz continuous gradients, we always have the following relation (see Theorem 2.1.5 in \cite{Nesterov-book}):
\[
F(u)\leq F(v)+\langle \nabla F(v), u-v\rangle+\frac{L}{2}\|v-u\|^2
\]
for any $u,v\in\mathbb{R}^d$.

Letting $u=\bar{x}^{k+1}$ and $v=\bar{x}^k$ in the preceding relation  yields
\begin{equation}\label{eq:F}
\begin{aligned}
&F(\bar{x}^{k+1})\\
&\leq F(\bar{x}^{k})+\langle \nabla F(\bar{x}^{k}), \bar{x}^{k+1}-\bar{x}^{k}\rangle +\frac{L}{2}\|\bar{x}^{k+1}-\bar{x}^{k}\|^2\\
&\leq F(\bar{x}^{k})+\langle \nabla F(\bar{x}^{k}), \gamma^k\bar\zeta_w^k-\gamma^k\bar{\xi}_w^k-\lambda^{k}\bar{g}^{k}\rangle \\
&\quad+\frac{L}{2}\|\gamma^k\bar\zeta_w^k-\gamma^k\bar{\xi}_w^k-\lambda^{k}\bar{g}^{k}\|^2,
\end{aligned}
\end{equation}
where in the second inequality we used the relation in (\ref{eq:bar_x}).

Subtracting $F(\theta^{\ast})$ on both sides of (\ref{eq:F})  and then taking the expectation (conditioned on $\mathcal{F}^k$) on both sides yield
\begin{equation}\label{eq:F_2}
\begin{aligned}
&\mathbb{E}\left[F(\bar{x}^{k+1})-F(\theta^{\ast})|\mathcal{F}^k\right]\\
&\leq F(\bar{x}^{k})-F(\theta^{\ast})-\langle \nabla F(\bar{x}^{k}),   \lambda^k{\bar g}^k\rangle \\
&\quad +\frac{L}{2} \mathbb{E}\left[\|\gamma^k\bar\zeta_w^k-\gamma^k\bar{\xi}_w^k-\lambda^{k}\bar{g}^{k}\|^2\right]\\
&\leq F(\bar{x}^{k})-F(\theta^{\ast})-\langle \nabla F(\bar{x}^{k}),   \lambda^k{\bar g}^k\rangle \\
&\quad+\frac{3L}{2} (\lambda^k)^2\|{\bar g}^k\|^2+ \frac{3L}{2}(\gamma^k)^2\mathbb{E}\left[\|\bar{\zeta}_w^k\|^2\right]\\
&\quad+\frac{3L}{2}(\gamma^k)^2\mathbb{E}\left[\|\bar{\xi}_w^k\|^2\right] \\
& \leq F(\bar{x}^{k})-F(\theta^{\ast})-\langle \nabla F(\bar{x}^{k}),   \lambda^k{\bar g}^k\rangle \\
&\quad+\frac{3L}{2} (\lambda^k)^2\|{\bar g}^k\|^2+ \frac{3L}{2}(\gamma^k)^2 \sum_{i,j} (R_{ij}\sigma_{\zeta,j}^{k})^2\\
&\quad+\frac{3L}{2}(\gamma^k)^2 \sum_{i,j} (C_{ij}\sigma_{\xi,j}^{k})^2.
\end{aligned}
\end{equation}
Next we bound the inner product term. Using the relationship $-\la a,b\ra=\frac{\|a-b\|^2-\|a\|^2-\|b\|^2}{2}$ valid for any vectors $a$ and $b$, one obtains
\begin{equation}\label{eq:inner_product}
\begin{aligned}
&-\langle \nabla F(\bar{x}^{k}),   \lambda^k{\bar g}^k\rangle\\
&=\frac{\lambda^k}{2}\left(\|\nabla F(\bar{x}^{k})-{\bar g}^k \|^2- \|\nabla F(\bar{x}^{k})\|^2- \|{\bar g}^k\|^2\right)\\
&\leq \frac{\lambda^k}{2}\left(\left\|\frac{1}{m}\sum_{i=1}^{m}(\nabla f_i(\bar{x}^{k})-\nabla f_i(x_i^k)  )\right\|^2- \|\nabla F(\bar{x}^{k})\|^2\right.\\
&\qquad \left.- \|{\bar g}^k\|^2\right)\\
&\leq \frac{\lambda^kL^2}{2m}\sum_{i=1}^{m} \| x_i^k-\bar{x}^{k} \|^2- \frac{\lambda^k}{2}\|\nabla F(\bar{x}^{k})\|^2- \frac{\lambda^k}{2}\|{\bar g}^k\|^2\\
&= \frac{\lambda^kL^2}{2m} \| {\bf x}^k-{\bf 1}\bar{x}^{k} \|_2^2- \frac{\lambda^k}{2}\|\nabla F(\bar{x}^{k})\|^2- \frac{\lambda^k}{2}\|{\bar g}^k\|^2.
\end{aligned}
\end{equation}

Plugging (\ref{eq:inner_product}) into (\ref{eq:F_2}) leads to
\begin{equation}\label{eq:F_final}
\begin{aligned}
&\mathbb{E}\left[F(\bar{x}^{k+1})-F(\theta^{\ast})|\mathcal{F}^k\right]\\
&\leq F(\bar{x}^{k})-F(\theta^{\ast}) +\frac{\lambda^kL^2}{2m}\| {\bf x}^k-{\bf 1}\bar{x}^{k} \|_2^2- \frac{\lambda^k}{2}\|\nabla F(\bar{x}^{k})\|^2\\
&\quad- \left(\frac{\lambda^k-3L (\lambda^k)^2}{2}\right)\|{\bar g}^k\|^2
 + \frac{3L}{2}(\gamma^k)^2 \sum_{i,j} (R_{ij}\sigma_{\zeta,j}^{k})^2\\
 &\quad + \frac{3L}{2}(\gamma^k)^2 \sum_{i,j} (C_{ij}\sigma_{\xi,j}^{k})^2.
\end{aligned}
\end{equation}

Step IV: We combine Steps I-III and prove the theorem.

Defining
\[
\bv^k=\big[F(\bar
x^{k+1})-F(\theta^{\ast}),\|{\bf x}^{k}-{\bf 1}\bar{x}^{k}\|_R^2,\|{\bf s}^{k}-v { \bar s}^{k}\|_C^2\big]^T,
\]
we
have the following relations from (\ref{eq:s_final}),
(\ref{eq:x-bar_x_final}), and (\ref{eq:F_final}):
\begin{equation}\label{eq:stacked}
\begin{aligned}
\mathbb{E}\left[\bv^{k+1}|\mathcal{F}^k\right]\leq(V^k+A^k)\bv^k-H^k\left[\begin{array}{c} \left\|\nabla
F(\bar x^k)\right\|^2\cr \|\bar g^k\|^2
\end{array}\right]+B^k,
\end{aligned}
\end{equation}
where
\[
\begin{aligned}
&V^k=\left[\begin{array}{ccc} 1 &
\frac{\delta^2_{2,R}\lambda^kL^2}{2m}&  0  \cr
0& 1-\gamma_1^k\rho_R&
\frac{2\gamma^k\|\Pi_UC\|_R^2\delta^2_{R,C}}{\rho_R}\cr 0&0
 &
1-\gamma_2^k\rho_c
\end{array}\right],\\
&A^k=\left[\begin{array}{ccc} 0 &
0&  0 \cr
a_1^k & a_2^k&0
\cr
 a_3^k&
a_4^k&
0
\end{array}\right],\\
&H^k=\left[\begin{array}{cc}   \frac{\lambda^k}{2} &
\frac{\lambda^k-3 L (\lambda^k)^2}{2}\cr 0&0\cr
0&0
\end{array}\right],\: B^k=\left[\begin{array}{c} b_1^k
\cr
b_2^k
\cr
b_3^k
\end{array}\right],
\end{aligned}
\]
with
\[
\begin{aligned}
a_1^k&=\frac{16mL(\lambda^k)^2\|\Pi_U\|_R^2\delta^2_{R,2}}{\gamma^k\rho_R},\\
a_2^k&=\frac{4(\lambda^k)^2 L^2\|\Pi_U\|_R^2\delta^2_{R,2}}{\gamma^k\rho_R},\\
a_3^k&=\frac{8mL(\lambda^{k})^2\delta^2_{C,2}\|\Pi_v\|^2_C}{\gamma^k\rho_c},\\
a_4^k&=\frac{2L^2(\lambda^{k})^2\delta^2_{C,2}\delta^2_{2,R}\|\Pi_v\|^2_C}{\gamma^k\rho_c},
\end{aligned}
\]
\[
\begin{aligned}
b_1^k&=\frac{3L}{2}(\gamma^k)^2 \sum_{i,j} (R_{ij}\sigma_{\zeta,j}^{k})^2 + \frac{3L}{2}(\gamma^k)^2 \sum_{i,j} (C_{ij}\sigma_{\xi,j}^{k})^2,\\
b_2^k&=\frac{8(\lambda^k)^2\|\Pi_U\|_R^2\delta^2_{R,2}}{\gamma^k\rho_R}   \sum_{i=1}^{m} \|\nabla f_i(\theta^\ast)\|^2    \\
 &\quad+ 2(\gamma^k)^2\|\Pi_u\|_R^2   \delta^2_{R,2} \sum_{i,j} (R_{ij}\sigma_{\zeta,j}^{k})^2   \\
 &\quad+2(\gamma^k)^2\|\Pi_U\|_R^2 \delta^2_{R,2} \sum_{i,j} (C_{ij}\sigma_{\xi,j}^{k})^2, \\
b_3^k&=\frac{4(\lambda^{k})^2\delta^2_{C,2}\|\Pi_v\|^2_C}{\gamma^k\rho_c} \sum_{i=1}^{m} \|\nabla f_i(\theta^\ast)\|^2\\
  &\quad+ (\gamma^k)^2\delta^2_{C,2}\|\Pi_v\|_C^2   \sum_{i,j} (C_{ij}\sigma_{\xi,j}^k)^2.
\end{aligned}
\]

Under Assumption \ref{assumption:dp-noises-intrack},  and the conditions that $(\gamma^k)^2$  $\frac{(\lambda^k)^2}{\gamma^k}$ are summable in the theorem statement,  it follows  that all entries of
the matrix $B^k$  are summable almost surely.  By  defining $\hat{b}^k$ as  the maximum element of $B^k$, we have $B^k\leq \hat{b}^k {\bf 1}$. Therefore,
  $\mathbb{E}\left[ F(\bar x^k)-F(\theta^{\ast})|\mathcal{F}^k\right]$,
$\mathbb{E}\left[\|{\bf x}^{k}-{\bf 1}\bar x^{k}\|_R^2|\mathcal{F}^k\right]$, and
$\mathbb{E}\left[\|{\bf s}^{k}-v{ \bar s}^{k}\|_C^2|\mathcal{F}^k\right]$ for the iterates generated by the proposed algorithm satisfy the conditions of Theorem \ref{Theorem:general_gradient_tracking} and,  hence, the results of Theorem \ref{Theorem:general_gradient_tracking} hold.
\end{proof}

{
\begin{Remark 1}
The requirement on the decaying-factor $\gamma^k$ and stepsize $\lambda^k$ in the statement  of Theorem \ref{th:deterministic} can be satisfied, for example, by setting $\gamma^k=\mathcal{O}(\frac{1}{k^a})$ and $\lambda^k=\mathcal{O}(\frac{1}{k^b})$ with $a,b\in\mathbb{R}$ satisfying $0.5<a<b\leq 1$ and $2b-a>1$. For example, setting $\gamma^k=\frac{c_1}{1+c_2k^{\iota}}$ and $\lambda^k=\frac{c_3}{1+c_4k}$ will satisfy the conditions for any   exponent $0.5<\iota<1$, and positive coefficients   $c_1$, $c_2$, $c_3$, and $c_4$.
\end{Remark 1}

\begin{Remark 1}
 Using the relationship in (\ref{eq:bar_y}) and the definitions of $\bar{s}^k$, $\bar{\xi}^k$, and $\bar{g}^k$ in (\ref{eq:definiton1}) and (\ref{eq:definition2}), one can obtain

\[
\bar{s}^{k+1}-\bar{s}^k=\gamma^k\bar{\xi}^k_w+\lambda^k\bar{g}^k,
\]
i.e.,  $\bar{s}^{k+1}-\bar{s}^k$  tracks the global gradient. Combined with the proven result in Theorem \ref{th:deterministic} that
 all $s_i^k$ converge to each other, and hence to the average $\bar{s}^k$ of all $s_i^k$, one can deduce  that $s_i^{k+1}-s_i^{k}$ in (\ref{eq:update}) of the proposed Algorithm 1  indeed  tracks the global gradient.

\end{Remark 1}

\section{Online Estimation of the left eigenvector}\label{se:estimation_vector}

In Algorithm 1, when the communication graph $\mathcal{G}_R$ is not balanced, a preprocessing approach can be used to estimate the left eigenvector $u^T$. In this section, inspired by the online eigenvector estimation algorithm in \cite{mai2016distributed}, we propose Algorithm 2 below, which allows individual agents to estimate the left eigenvector $u^T$ locally on the fly while  updating their optimization iterations in a distributed manner:

\noindent\rule{0.49\textwidth}{0.5pt}
\noindent\textbf{Algorithm 2: Robust gradient-tracking
based distributed optimization with eigenvector estimation}

\vspace{-0.2cm}\noindent\rule{0.49\textwidth}{0.5pt}
\begin{enumerate}[wide, labelwidth=!, labelindent=0pt]
    \item[] Parameters: Stepsize $\lambda^k$ and a decaying factor $\gamma^k$  to suppress information-sharing noise;
    %,  in-bound mixing/pulling weighs $R_{ij}>0$ for all $j\in\mathbb{N}_{R,i}^{\rm in}$, and out-bound pushing weights $C_{l,i}>0$ for all $l\in\mathbb{N}_{C,i}^{\rm out}$, otherwise $R_{ij}=C_{li}=0$;
    \item[] Every agent $i$ maintains two states  $x_i^k$ and
    $s_i^k$, which are initialized randomly with $x_i^0\in\mathbb{R}^d$ and $s_i^0\in\mathbb{R}^d$. Every agent $i$ also maintains an eigenvector-estimation parameter $z_i^k\in\mathbb{R}^m$ initialized with $z_i^0={\bf e}_i\in\mathbb{R}^m$ where ${\bf e}_i$ has the $i$th element equal to one and all other elements equal to zero.
    \item[] {\bf for  $k=1,2,\cdots$ do}
    \begin{enumerate}
        %\item Every agent $j$ injects zero-mean DP-noises $\zeta_j^k$ and $\xi_j^k$  to its states      $x_j^k$ and $y_j^k$, respectively.
        \item Agent $i$ pushes $s_i^k$ to each agent
        $l\in\mathbb{N}_{C,i}^{\rm out}$, which will be received as $s_i^k+\xi_i^k$ due to   information-sharing noise.  And agent $i$ pulls $x_j^k$ from each $j\in\mathbb{N}_{R,i}^{\rm in}$, which will be received as $x_j^k+\zeta_j^k$ due to information-sharing noise. Here the subscript  $R$ or $C$ in neighbor sets indicates the neighbors with respect to  the graphs induced by these matrices. Agent $i$ also pulls $z_j^k$ from each $j\in\mathbb{N}_{R,i}^{\rm in}$.
         \item Agent $i$ chooses $\gamma^k>0$  satisfying
        $1+\gamma^kR_{ii}>0$ and $1+\gamma^kC_{ii}>0$ with  $R_{ii}$ and $C_{ii}$ defined in (\ref{eq:diagonal_entries}).
        Then, agent $i$ updates its states as follows:
        \begin{equation}\label{eq:update_2}
        \begin{aligned}
s_i^{k+1}=&(1+\gamma^kC_{ii})s_i^k+\gamma^k\sum_{j\in \mathbb{N}_{C,i}^{\rm
            in}}C_{ij}(s_j^k+\xi_j^k)\\
            &\quad + \lambda^k\nabla f_i(x_i^{k}),\\
x_i^{k+1}=&(1+\gamma^k R_{ii})x_i^k+\gamma^k\sum_{j\in \mathbb{N}_{R,i}^{\rm
            in}}R_{ij}(x_j^k+\zeta_j^k)\\
            &\quad-\frac{s_i^{k+1}-s_i^{k}}{mz^k_{ii}},\\
z_{i}^{k+1}=&z_i^k+ \sum_{j\in \mathbb{N}_{R,i}^{\rm
            in}}R_{ij}(z_j^k-z_i^k),
        \end{aligned}
         \end{equation}
         where $z^k_{ii}$ denotes the $i$th element of $z_i^k$.

                  \item {\bf end}
    \end{enumerate}
\end{enumerate}
\vspace{-0.1cm} \rule{0.49\textwidth}{0.5pt}

In Algorithm 2, every agent uses the third update in (\ref{eq:update_2}) to locally estimate the left eigenvector of $I+\gamma^kR$. (Note that as discussed in Remark \ref{re:eigenvector_time_invariant}, the left eigenvector of $I+\gamma^kR$ is time-invariant and independent of $\gamma^k$. Also note that the update obtains an estimated eigenvector with row sum equal to one, and thus  we scale the estimate by $m$ to obtain $u^T$ whose row sum is required to be $m$.) Therefore, every agent $i$ can use its local estimate $z_i^k$ of the left eigenvector, which avoids using global information of $u^T$ in Algorithm 1. It is worth noting that since  $z_i^k$ does not contain sensitive information, there is no need to add information-sharing noise to them to enable differential privacy. Moreover, the dimension of $z_i^k$ is equal to the size of the network $m$. Thus,  even in the case where the communication channel is noisy or coarse quantization is used for $s$-iterates and $x$-iterates, special effort (e.g., error-correction coding \cite{clark2013error} or high-precision quantization) can be exploited to ensure that  shared $z_i^k$ messages are not   contaminated by noises. Note that such special effort may not  be feasible for the sharing of optimization variables ($s$-iterates and $x$-iterates) since the dimension of optimization variables  can scale up to  hundreds of millions in deep learning applications  \cite{huang2017densely}, which makes the cost for error-correction coding or high-precision quantization  prohibitively high.

Next, we prove that Algorithm 2 can still ensure almost sure convergence of all agents to an optimal solution. To this end, we first characterize the estimation error of the eigenvector estimator:
\begin{Lemma 7}\label{Le:geometric_convergence}
Under Assumption \ref{Assumption:push_pull topology},  the iterates $z_i^k$   in (\ref{eq:update_2}), after scaled by $m$, converge to the  left eigenvector $u^T=\left[u_1,\,u_2,\,\cdots,u_m\right]^T$ of $I+\gamma^k R$ with a geometric rate, i.e., there exist $C>0$ and $p\in(0,\,1)$ satisfying the following inequality for any $i\in[m]$ and $k\geq 0$:
\begin{equation}\label{eq:error_geometric}
\left|\frac{1}{mz_{ii}^k}-\frac{1}{u_i}\right|\leq Cp^k,
\end{equation}
where $z_{ii}^k$ denotes the $i$th element of $z_i^k$.
\end{Lemma 7}
\begin{proof}
From \cite{mai2016distributed}, we know that there exist $C_1>0$ and $p\in(0,\,1)$ such that $\left|mz_{ii}^k-u_i\right|\leq C_1p^k$ holds under the given conditions. According to \cite{mai2016distributed}, we also know that $u_i$ and $z^k_{ii}$ are strictly positive numbers. Therefore, using the relation $\left|\frac{1}{mz_{ii}^k}-\frac{1}{u_i}\right|=\frac{\left|mz_{ii}^k-u_i\right|}{mz^k_{ii}u_i}$, we know that  there exist $C>0$ such that (\ref{eq:error_geometric}) holds.
\end{proof}

Based on Lemma \ref{Le:geometric_convergence}, we can prove the almost sure convergence of all agents to an optimal solution following the line of reasoning of Theorem \ref{th:deterministic}:

\begin{Theorem 5}\label{th:with_estimation}
 Let  Assumption 1, Assumption~\ref{Assumption:push_pull topology}, and
Assumption~\ref{assumption:dp-noises-intrack} hold.
If  $\{\gamma^k\}$    and $\{\lambda^k\}$ satisfy
$
 \sum_{k=0}^\infty \gamma^k=\infty$,  $\sum_{k=0}^\infty (\gamma^k)^2<\infty$, $\sum_{k=0}^\infty \lambda^k=\infty$, $\sum_{k=0}^\infty \frac{(\lambda^k)^2}{\gamma^k}<\infty$, and $\lim_{k\to\infty}\lambda^k/\gamma^k=0$,
then the results of Theorem~\ref{Theorem:general_gradient_tracking}
hold for  Algorithm 2.
\end{Theorem 5}
\begin{proof}
 The proof follows the derivation of Theorem \ref{th:deterministic}. Since the eigenvalue estimation process does not affect  the dynamics of $s_i^k$,  the relation   for $\|{\bf s}^k-v{ \bar{s}}^k \|_C$ in Step I of Theorem \ref{th:deterministic} still holds for Algorithm 2. Thus, we only need to show that we can establish relations for
 $\|{\bf x}^k-{\bf 1}{ \bar{x}}^k \|_R$ and $\mathbb{E}\left[ F(\bar x^k)-F(\theta^{\ast})|\mathcal{F}^k\right]$ that are  similar to those in Theorem \ref{th:deterministic}.

 %Relationship for $\|{\bf x}^k-{\bf 1}{ \bar{x}}^k \|_R$.

Similar to (\ref{eq:bar_x}), denoting $U$ as ${\rm diag}(u_1,\,u_2,\,\cdots,u_m)$ and $Z^k$ as ${\rm diag}(mz^k_{11},\,mz^k_{22},\,\cdots,mz^k_{mm})$, with $z_{ii}^k$ denoting the $i$th element of $z_i^k$, we can obtain the following relationship for the $x$-iterates in Algorithm 2:
\begin{equation}\label{eq:bar_x_2}
\begin{aligned}
&\bar{x}^{k+1}\\
=& \frac{u^T}{m}{\bf x}^{k+1}= \frac{u^T}{m}( R^k  {\bf x}^k+\gamma^k {\pmb \zeta}_w^k- (Z^k)^{-1} ({\bf s}^{k+1}-{\bf s}^{k}))\\
=&\frac{u^T}{m}\left( R^k  {\bf x}^k+\gamma^k {\pmb \zeta}_w^k- \hspace{-0.1cm}\left(U^{-1}\hspace{-0.1cm}+\hspace{-0.1cm}(Z^k)^{-1}\hspace{-0.1cm}-\hspace{-0.1cm}U^{-1}\right) ({\bf s}^{k+1}-{\bf s}^{k})\right)\\
=&\bar{x}^k+\gamma^k\bar\zeta_w^k-\frac{{\bf 1}^T}{m} ({\bf s}^{k+1}-{\bf s}^{k})\\
&-\frac{{u}^T((Z^k)^{-1}-U^{-1})}{m}({\bf s}^{k+1}-{\bf s}^{k})).
%&=\bar{x}^k+\gamma^k\bar\zeta_w^k-({\bar s}^{k+1}-{\bar s}^{k})-(u_{e}^k)^T({\bf s}^{k+1}-{\bf s}^{k})\\
%&=\bar{x}^k+\gamma^k\bar\zeta_w^k-\gamma^k\bar{\xi}_w^k-\lambda^{k}\bar{g}^{k}- \gamma^k(u_{e}^k)^TC({\bf s}^k-v\bar{s}^k)\\
%&\qquad  +\gamma^k(u_{e}^k)^T {\pmb\xi}_w^k+\lambda^k(u_{e}^k)^T{\bf g}^{k}
\end{aligned}
\end{equation}

Hence, following the line of reasoning in the proof of Theorem \ref{th:deterministic}, we can obtain
\begin{equation}\label{eq:x_bar_xalgorithm2}
\begin{aligned}
{\bf x}^{k+1}-{\bf 1}\bar{x}^{k+1}=&\bar{R}^k ({\bf x}^{k}-{\bf 1}\bar{x}^{k})-\gamma^k\Pi_U C({\bf s}^k-v\bar{s}^k)\\
&  -\gamma^k\Pi_U {\pmb\xi}_w^k  -\lambda^k\Pi_U{\bf g}^{k}  +\gamma^k\Pi_u{\pmb \zeta}_w^k\\
& +\gamma^k\Pi_U^e({\bf s}^k-v\bar{s}^k)+\gamma^k\Pi_U^e {\pmb\xi}_w^k +\lambda^k\Pi_U^e{\bf g}^{k},
\end{aligned}
\end{equation}
where $\Pi_U^e=(I-\frac{{\bf 1}u^T}{m})((Z^k)^{-1}-U^{-1})$ and we have used the relationship in (\ref{eq:s_s}) in the last equality.

In (\ref{eq:x_bar_xalgorithm2}), the last three terms on the right hand side    correspond  to the influence of introducing the eigenvector estimator. Given that the elements of $(Z^k)^{-1}-U^{-1}$ diminish  with a geometric rate according to Lemma \ref{Le:geometric_convergence}, we deduce  that the coefficient sequences for these terms are all summable, and hence they will only introduce terms with summable  coefficients in the relationship for  $\|{\bf x}^k-{\bf 1}{ \bar{x}}^k \|_R$, which will not affect the almost sure convergence. The same reasoning applies to the   dynamics of $\bar{x}^{k+1}-\bar{x}^k$. More specifically,  compared with Algorithm 1, Algorithm 2's eigenvector estimator (the last item on the right hand side of (\ref{eq:bar_x_2})) introduces three extra terms $\gamma^k\frac{{u}^T((Z^k)^{-1}-U^{-1})}{m}C({\bf s}^k-v\bar{s}^k)$, $\gamma^k\frac{{u}^T((Z^k)^{-1}-U^{-1})}{m} {\pmb\xi}_w^k$, and $\lambda^k\frac{{u}^T((Z^k)^{-1}-U^{-1})}{m}{\bf g}^{k}$ according to (\ref{eq:s_s}). From Lemma \ref{Le:geometric_convergence}, we know that  their coefficients  all decrease  with a geometric rate and hence are all summable. Therefore, these three extra terms only introduce items that have summable coefficient sequences in the relationship for $F(\bar{x}^{k})-F(\theta^{\ast})$, which will not affect the almost sure convergence.

In summary, we have that introducing the eigenvector estimator adds terms with summable coefficient sequences in the final   inequality in (\ref{eq:stacked}), and hence will not affect the almost sure convergence results in Theorem \ref{th:deterministic}. Therefore, we can still prove that the   iterates generated by Algorithm 2 satisfy the conditions of Theorem \ref{Theorem:general_gradient_tracking} and,  hence, the results of Theorem \ref{Theorem:general_gradient_tracking} hold for Algorithm 2.
\end{proof}
}

\section{Extension to Distributed Stochastic Gradient Methods}\label{se:SGD}
In many distributed optimization applications, individual agents do not have access to
the precise gradient and hence have to use  {\it noisy}
local gradients for optimization.
For example, in modern machine learning on massive
datasets,  evaluating the
precise gradient using all available data can be extremely expensive
in computation or even practically infeasible. So individual agents usually only compute   inexact  estimates  of the true gradients using a portion of the
data points available to them \cite{koloskova2019decentralized}.
  Furthermore, in the
era of Internet of Things, which connect  massive low-cost sensing
and communication devices, the data fed to optimization computations
are usually subject to measurement noises
\cite{xin2020decentralized}. In this  section, we prove that the proposed algorithm can ensure all agents' almost sure convergence to an optimal solution even when the gradients are noisy.

As in most existing results on stochastic gradient methods, we make the following standard assumption on the stochasticity of individual agents' local gradients:
\begin{Assumption 4}
Every individual agent's local gradient  $g_i^k$ is an unbiased estimate of the true gradient $\nabla f_i(x_i^k)$ and has bounded variance, i.e.,
\[
\mathbb{E}\left[g_i^k\right]=\nabla
f_i(x_i^k),\:\forall i
\]
\[
 \mathbb{E}\left[\|g_i^k-\nabla
  f_i(x)\|^2\right]\leq \sigma^2,\: \forall i, x
\]
where $\sigma$ is some  positive constant.
\end{Assumption 4}

\begin{Theorem 1}
 Let  Assumptions 1-4 hold.
If  $\{\gamma^k\}$    and $\{\lambda^k\}$ satisfy
$
 \sum_{k=0}^\infty \gamma^k=\infty$,  $\sum_{k=0}^\infty (\gamma^k)^2<\infty$, $\sum_{k=0}^\infty \lambda^k=\infty$, $\sum_{k=0}^\infty \frac{(\lambda^k)^2}{\gamma^k}<\infty$, and $\lim_{k\to\infty}\lambda^k/\gamma^k=0$,
then the results of Theorem~\ref{Theorem:general_gradient_tracking}
hold for the proposed Algorithm 1 and Algorithm 2  even when individual agents have access to  only stochastic estimates of their true gradients.
\end{Theorem 1}
\begin{proof}
{
We use Algorithm 1 as an example to prove the results. Similar derivations apply to Algorithm 2 as well.

The  goal is still to establish the relationship in~(\ref{eq-fine}), with the $\sigma$-field
$\mathcal{F}^k =\{x_i^\ell,s_i^\ell; 0\le \ell\le k,\ i\in[m]\}$.
To this end, we organize  the derivations into four
steps: in Step I, Step II, and Step III, we establish respectively  relations
for $\|{\bf s}^k-v{ \bar{s}}^k \|_C$, $\|{\bf x}^k-{\bf 1}{ \bar{x}}^k \|_R$, and $\mathbb{E}\left[ F(\bar x^k)-F(\theta^{\ast})|\mathcal{F}^k\right]$ for the iterates generated by the proposed algorithm. In Step IV, we use them to
show that~(\ref{eq-fine}) of Theorem~\ref{Theorem:general_gradient_tracking} holds.

Step I: Relationship for $\|{\bf s}^k-v{ \bar{s}}^k \|_C$.

Since  the noise on gradients can be grouped into the noise term $\xi_i^k$, following the same procedure as in Theorem \ref{th:deterministic}, we can obtain a relation similar to (\ref{eq:s_final}):
\begin{equation}\label{eq:s_final_ext}
\begin{aligned}
& \mathbb{E}\left[ \left\|{\bf s}^{k+1}-v{ \bar{s}}^{k+1}\right\|^2_C|\mathcal{F}^k\right] \leq(1-\gamma^k\rho_c)\left\| {\bf s}^k -v{ \bar{s}}^{k} \right\|_C^2 \\
&\qquad+\frac{2L^2(\lambda^{k})^2\delta^2_{C,2}\|\Pi_v\|^2_C}{\gamma^k\rho_c}  \|{\bf x}^k- {\bf 1}\bar{x}^k \|_2^2\\
&\qquad+\frac{8mL(\lambda^{k})^2\delta^2_{C,2}\|\Pi_v\|^2_C}{\gamma^k\rho_c} ( F(\bar{x}^k)-F(\theta^\ast))   \\
&\qquad+\frac{4(\lambda^{k})^2\delta^2_{C,2}\|\Pi_v\|^2_C}{\gamma^k\rho_c} \sum_{i=1}^{m} \|\nabla f_i(\theta^\ast)\|^2\\
  &\qquad+ (\gamma^k)^2\delta^2_{C,2}\|\Pi_v\|_C^2   \sum_{i,j} (C_{ij}\sigma_{\xi,j}^k)^2\\
  &\qquad+ m^2(\lambda^k)^2\delta^2_{C,2}\|\Pi_v\|_C^2 \sigma^2,
\end{aligned}
\end{equation}
where the last term corresponds to the influence caused by the stochasticity in local gradients.

Step II: Relationship for $\|{\bf x}^k-{\bf 1}{ \bar{x}}^k \|_R$.

Still following the derivations in Theorem \ref{th:deterministic}, we have that the stochasticity in local gradients will affect the term $\|{\bf g}^{k}\|^2_R$ in (\ref{eq:Algorithm1_x_error}). More specifically, after taking conditional expectation, $\|{\bf g}^{k}\|^2_R$ will become $\|{\bf g}^{k}\|^2_R+m^2\delta_{R,2}^2\sigma^2$, and hence (\ref{eq:x-bar_x_final})  becomes
\begin{equation}\label{eq:x-bar_x_final_ext}
 \begin{aligned}
 \begin{aligned}
& \mathbb{E}\left[ \|{\bf x}^{k+1}-{\bf 1}\bar{x}^{k+1}\|_R^2 |\mathcal{F}^k\right]\\
&\leq \left(1-\gamma^k\rho_R+\frac{4(\lambda^k)^2 L^2\|\Pi_U\|_R^2\delta^2_{R,2}}{\gamma^k\rho_R}\right) \|{\bf x}^{k}-{\bf 1}\bar{x}^{k}\|_R^2\\
&\quad+\frac{2\gamma^k\|\Pi_UC\|_R^2}{\rho_R} \|{\bf s}^k-v\bar{s}^k\|_R^2 \\
 &\quad+
 \frac{16mL(\lambda^k)^2\|\Pi_U\|_R^2\delta^2_{R,2}}{\gamma^k\rho_R} (     F(\bar{x}^k)-F(\theta^\ast)) \\
 &\quad+
 \frac{8(\lambda^k)^2\|\Pi_U\|_R^2\delta^2_{R,2}}{\gamma^k\rho_R}   \sum_{i=1}^{m} \|\nabla f_i(\theta^\ast)\|^2    \\
 &\quad+ \frac{2m^2(\lambda^k)^2\|\Pi_U\|_R^2\delta^2_{R,2}}{\gamma^k\rho_R}  \sigma^2     \\
 &\quad+ 2(\gamma^k)^2\|\Pi_u\|_R^2   \delta^2_{R,2} \sum_{i,j} (R_{ij}\sigma_{\zeta,j}^{k})^2   \\
 &\quad+2(\gamma^k)^2\|\Pi_U\|_R^2 \delta^2_{R,2} \sum_{i,j} (C_{ij}\sigma_{\xi,j}^{k})^2.
 \end{aligned}
 \end{aligned}
\end{equation}

Step III: Relationship for $F(\bar{x}^{k})-F(\theta^{\ast})$.

Following the derivations in Theorem \ref{th:deterministic}, we have that the stochasticity in local gradients affects $\bar{g}^k$, which will be subject to noise with variance $\sigma^2$. More specifically, we have that (\ref{eq:F_final}) becomes
\begin{equation}\label{eq:F_final_ext}
\begin{aligned}
&\mathbb{E}\left[F(\bar{x}^{k+1})-F(\theta^{\ast})|\mathcal{F}^k\right]\\
&\leq F(\bar{x}^{k})-F(\theta^{\ast}) +\frac{\lambda^kL^2}{2m}\| {\bf x}^k-{\bf 1}\bar{x}^{k} \|_2^2- \frac{\lambda^k}{2}\|\nabla F(\bar{x}^{k})\|^2\\
&\quad- \left(\frac{\lambda^k-3L (\lambda^k)^2}{2}\right)\|{\bar g}^k\|^2
 +  \frac{3L}{2}(\lambda^k)^2\sigma^2 \\
 &\quad + \frac{3L}{2}(\gamma^k)^2 \sum_{i,j} (R_{ij}\sigma_{\zeta,j}^{k})^2+ \frac{3L}{2}(\gamma^k)^2 \sum_{i,j} (C_{ij}\sigma_{\xi,j}^{k})^2.
\end{aligned}
\end{equation}

Step IV: We combine Steps I-III and prove the theorem, which involves arguments  exactly the same as in the proof of Theorem \ref{th:deterministic}. In fact,  following the derivation in Theorem \ref{th:deterministic}, we can obtain  that the stochasticity of gradients will only affect the matrix $B^k$ in (\ref{eq:stacked}), which will still be summable. Therefore, we can arrive at the same conclusion as in Theorem \ref{th:deterministic} even when local gradients are stochastic.
}

% Following the derivations in the proof of Theorem \ref{th:deterministic}, it can be seen that in general, the noise on gradient can be grouped into the noise term $\xi_i^k$. There is only one exception in the relationship of $\|{\bf g}^k\|^2$ in (\ref{eq:g_norm}). In fact, the stochasticity of $g_i^k$ will bring a new constant term $\sigma^2$. However, given that 1) $\sigma^2$ is bounded; and 2) $\|{\bf g}^k\|^2$ is factored into the final relationship after multiplying by $\frac{(\lambda^k)^2}{\gamma^k}$ which is summable (see (\ref{eq:s_1}) and (\ref{eq:x-bar_x_2})), we know that the final influence of gradient noise on convergence iterates is still a summable term. In summary, we can obtain  that the stochasticity of gradients will only affect the matrix $B^k$ in (\ref{eq:stacked}), which will still be summable. Therefore, we can arrive at the same conclusion as in Theorem \ref{th:deterministic} even when local gradients are stochastic.
\end{proof}

\section{Numerical Simulations}\label{se:simulations}

In this section, we evaluate the performance of the proposed  distributed optimization algorithm  within the context of  a distributed estimation problem.

We consider a canonical distributed estimation problem where a
 network of $m$ sensors collectively estimate an  unknown parameter
$\theta\in\mathbb{R}^d$. More specifically, we assume that each
sensor $i$ has a noisy measurement  of the parameter,
$z_{i}=M_i\theta+w_{i}$, where $M_i\in\mathbb{R}^{s\times d}$ is the
measurement matrix of agent $i$ and $w_{i}$ is Gaussian measurement noise of unit variance.
Then the maximum likelihood  estimation of parameter $\theta$  can be solved using the
optimization problem formulated as
(\ref{eq:optimization_formulation1}), with each $f_i(\theta)$ given
as
$
f_i(\theta)=\|z_i-M_i\theta\|^2+\varsigma\|\theta\|^2
$,
where $\varsigma$ is a  regularization parameter \cite{xu2017convergence}.

\begin{figure}
\includegraphics[width=0.53\textwidth]{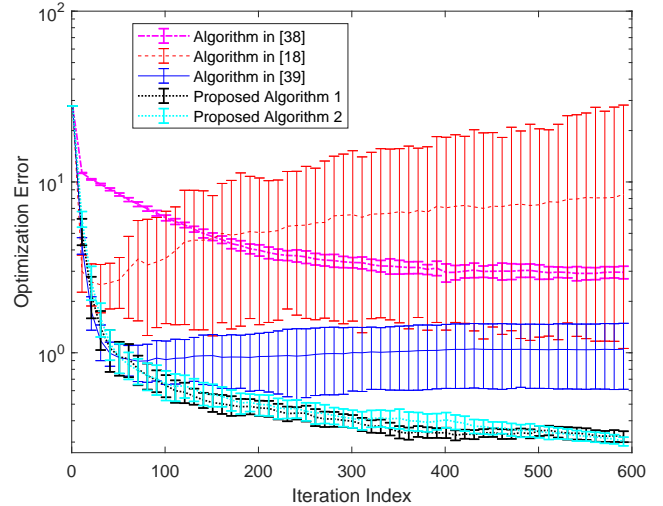}
    \caption{Comparison of the proposed algorithms  with the conventional Push-Pull algorithm in \cite{pu2020push}, the algorithm in \cite{pu2020robust} that can avoid noise accumulation when coupling matrices are constant, and the algorithm in \cite{wang2022tailoring} that combines conventional Push-Pull with decaying factors.}
    \label{fig:comparison}
\end{figure}

In the numerical experiments, we set the number of agents (sensors) to {$m=100$ and adopt  a random interaction graph. To ensure that the random interaction graph is strongly connected, we first arrange the 100 agents on a   ring and then add
  a directed link  between any two nonadjacent nodes with probability 0.3.} In the evaluation, we set $s=3$ and $d=2$. To evaluate the performance of the proposed algorithms,  we inject  Gaussian based information-sharing noise $\zeta_i^k$    and $\xi_i^k$  on all shared $x_i^k$ and $s_i^k$, respectively.  Both $\zeta_i^k$ and $\xi_i^k$ have mean  0 and standard deviation  $\sigma_{\zeta,i}^k=\sigma_{\xi,i}^k=0.8$. We set the stepsize $\lambda^k$ and diminishing sequence $\gamma^k$ as $\lambda^k=\frac{0.02}{1+0.1k}$ and $\gamma^k=\frac{1}{1+0.1k^{0.6}}$, respectively, which satisfy the conditions in Theorem 2. We run our Algorithm 1 and Algorithm 2 for 100 times and calculate  the average as well as the variance of the optimization error {$\sum_{i=1}^{m}\|x_i^k-\theta^{\ast}\|$} as a function of the iteration index $k$.  The result for Algorithm 1 is given by the black  curve  and error bars in Fig. \ref{fig:comparison}, and the result for Algorithm 2 is given by the  cyan curve  and error bars in Fig. \ref{fig:comparison}. For comparison, we also run the conventional  Push-Pull algorithm in \cite{pu2020push} (which uses a constant stepsize and no decaying factor), the robust gradient-tracking algorithm proposed in \cite{pu2020robust} (which uses a constant stepsize and can avoid noise accumulation under constant inter-agent coupling),  and our recent result in \cite{wang2022tailoring} (which combines the conventional Push-Pull with  decaying factors). The stepsize for the conventional Push-Pull method in \cite{pu2020push} and the algorithm in \cite{pu2020robust} is set to a constant value $\lambda^k=0.02$, and the decaying factor and stepsize for \cite{wang2022tailoring} is set  the same as ours (note that \cite{wang2022tailoring} uses two decaying factors, and we set one of them equal to our decaying factor and the other one is selected according to the requirement therein). For all these three algorithms, we run the experiments for 100 times under the same information-sharing noise.    The evolution of the average optimization errors and variances for the three algorithms  are depicted by  the  curves  and error bars in orange, magenta, and blue, respectively, in Fig. \ref{fig:comparison}. It is clear that the proposed algorithms have   both   faster convergence speeds and better optimization accuracies compared with existing results. Furthermore, it can be seen that the variance of the optimization error for the conventional Push-Pull algorithm in \cite{pu2020push} indeed grows with time, which corroborates the accumulation of information-sharing noise in conventional gradient-tracking based algorithms. {It is also worth noting that for the approach in  \cite{pu2020robust} with a constant stepsize, although the theoretical analysis therein establishes that the expected optimization error converges linearly to a steady-state value, the actual optimization error may decrease with a slower rate. }

%\begin{figure}
%\center
%\includegraphics[width=0.22\textwidth]{graph_directed}
%    \caption{The interaction topology of the network}
%    \label{fig:topology}
%\end{figure}

\section{Conclusions}\label{se:conclusions}

The robustness of distributed optimization algorithms against information-sharing noise is becoming increasingly important due to the prevalence of channel noise, the existence  of quantization errors, and the demand for data perturbation/randomization for privacy protection. However,  gradient-tracking based distributed optimization, which is gaining increased traction due to its applicability to general directed graphs and fast convergence speed, is vulnerable to information-sharing noise. In fact,  in existing algorithms, information-sharing noise accumulates on  the global gradient estimate and its variance will even grow to infinity when the noise is persistent. We have proposed a new gradient-tracking based approach which can avoid information-sharing noise from accumulating in the global-gradient estimate. The approach is applicable even when   { the inter-agent interaction is} time-varying, enabling the incorporation of  a decaying factor  to gradually eliminate the influence of information-sharing noise, even when the noise is persistent. We have proved that with an appropriately chosen decaying factor, the proposed approach can guarantee all agents' almost sure convergence to an optimal solution for general convex objective functions with Lipschitz gradients, even in the presence of persistent information-sharing noise. The approach is also applicable when local gradients are subject to bounded noises as well, which is common in machine learning applications.  Numerical simulation results confirm that in the presence of information-sharing noise, the proposed approach
has  better optimization accuracy compared with  existing
counterparts.

{ We should note that a limitation of our approach is that it assumes time-invariant coupling topology. We plan to explore relaxation of this assumption  in future work. Moreover, in future work, we also plan to study whether  decaying  factors  can be incorporated into non-gradient based distributed optimization algorithms to  enable robustness against information-sharing noise. }

\bibliographystyle{IEEEtran}

\bibliography{reference1}

\end{document}